\documentclass[a4paper,11pt,reqno]{amsart}
\usepackage[a4paper,margin=0.9in]{geometry}
\usepackage{amsmath,amssymb,amsthm,bm,mathrsfs}
\usepackage{microtype,enumitem}
\usepackage[hidelinks]{hyperref}
\hypersetup{
 pdftitle={Critical Resonance and the Radial Branch in Maz'ya's Weighted Hessian--Gradient Inequality},
 pdfauthor={Suman Das},
 pdfsubject={Sharp critical asymptotics, exterior radiality, and symmetry breaking for Maz'ya's planar weighted Hessian--gradient inequality},
 pdfkeywords={weighted Hessian--gradient inequality, sharp constants, critical resonance, radial symmetry, symmetry breaking, Hardy--Leray inequality, Emden--Fowler transform}
}
\allowdisplaybreaks
\numberwithin{equation}{section}
\setlist{nosep,leftmargin=1.4em}
\usepackage{aliascnt}
\usepackage[nameinlink,noabbrev]{cleveref}
\newtheorem{theorem}{Theorem}[section]
\newaliascnt{proposition}{theorem}
\newtheorem{proposition}[proposition]{Proposition}
\aliascntresetthe{proposition}
\newaliascnt{lemma}{theorem}
\newtheorem{lemma}[lemma]{Lemma}
\aliascntresetthe{lemma}
\newaliascnt{corollary}{theorem}
\newtheorem{corollary}[corollary]{Corollary}
\aliascntresetthe{corollary}
\theoremstyle{remark}
\newaliascnt{remark}{theorem}
\newtheorem{remark}[remark]{Remark}
\aliascntresetthe{remark}
\crefname{theorem}{Theorem}{Theorems}
\Crefname{theorem}{Theorem}{Theorems}
\crefname{proposition}{Proposition}{Propositions}
\Crefname{proposition}{Proposition}{Propositions}
\crefname{lemma}{Lemma}{Lemmas}
\Crefname{lemma}{Lemma}{Lemmas}
\crefname{corollary}{Corollary}{Corollaries}
\Crefname{corollary}{Corollary}{Corollaries}
\crefname{remark}{Remark}{Remarks}
\Crefname{remark}{Remark}{Remarks}
\crefname{section}{Section}{Sections}
\Crefname{section}{Section}{Sections}
\newcommand{\R}{\mathbb R}
\newcommand{\T}{\mathbb T}
\newcommand{\cA}{\mathcal A}
\newcommand{\Frob}{\mathrm F}
\newcommand{\rad}{\mathrm{rad}}
\newcommand{\ang}{\mathrm{ang}}
\newcommand{\sgn}{\operatorname{sgn}}
\newcommand{\norm}[1]{\left\lVert#1\right\rVert}
\newcommand{\dd}{\,\mathrm d}
\title[Maz'ya's Weighted Hessian--Gradient Inequality]{Critical Resonance and the Radial Branch in\\Maz'ya's Weighted Hessian--Gradient Inequality}
\author{Suman Das}
\address{Department of Mathematics, Shantou University, Shantou, Guangdong 515063, P.~R.~China}
\email{suman@stu.edu.cn}
\date{}
\begin{document}
\begin{abstract}
We study the sharp constant $C(p,\gamma)$ in Maz'ya's planar weighted inequality
\[
 \int_{\R^2}|\nabla\psi|^p|x|^{p(\gamma-1)}\dd x
 \le C(p,\gamma)\int_{\R^2}|D^2\psi|_{\Frob}^p|x|^{p\gamma}\dd x,
\]
where $|\cdot|_{\Frob}$ denotes the Frobenius norm. For every fixed $1\le p<\infty$, with $a=\gamma-1+2/p$, we prove
\[
 C(p,\gamma)^{-1}=2^{p/2}|a|^p+O_p(|a|^{p+1})\qquad(a\to0).
\]
The leading term comes from the first angular harmonics after the Emden--Fowler transform. We also determine explicit exterior parameter rays on which the unrestricted constant equals the radial value $(1+a^2)^{-p/2}$. The proof uses an entropy identity and a homogeneous calibration. In the remaining transition region, we obtain symmetry-breaking intervals and a necessary and sufficient condition for nonnegativity of the periodic second variation at the radial state. An explicit example shows that global symmetry breaking can occur even when this second variation is nonnegative.
\end{abstract}
\maketitle

\section{Introduction and main results}\label{sec:intro}
Let $1\le p<\infty$ and $\gamma\in\R$. For a real-valued function $\psi\in C_c^\infty(\R^2)$, we denote
\[
 I_{p,\gamma}[\psi]=\int_{\R^2}|\nabla\psi(x)|^p|x|^{p(\gamma-1)}\dd x,
 \qquad
 J_{p,\gamma}[\psi]=\int_{\R^2}|D^2\psi(x)|_{\Frob}^p|x|^{p\gamma}\dd x.
\]
Here $|D^2\psi|_{\Frob}^2=\sum_{i,j=1}^2\psi_{x_ix_j}^2$ is the Frobenius norm. We consider the finite-energy class
\[
 \cA_{p,\gamma}=\left\{\psi\in C_c^\infty(\R^2):
 0<I_{p,\gamma}[\psi]<\infty,\quad J_{p,\gamma}[\psi]<\infty,\quad
 \nabla\psi(0)=0\ \text{if }\gamma<1-\frac2p\right\}.
\]
The best constant and its reciprocal are
\begin{equation}\label{eq:constant}
 C(p,\gamma)=\sup_{\psi\in\cA_{p,\gamma}}
 \frac{I_{p,\gamma}[\psi]}{J_{p,\gamma}[\psi]},
 \qquad
 \Lambda_p(a)=C(p,\gamma)^{-1},
 \qquad a=\gamma-1+\frac2p.
\end{equation}
We use the convention $\infty^{-1}=0$. Thus $\Lambda_p(a)$ is the infimum of the reversed quotient. Its numerator is positive for every admissible test: a smooth compactly supported function with zero Hessian is zero. Throughout the paper, we write
\[
 \beta=1-a,\qquad B=1+a^2.
\]

For $a\ne0$, the determination of \eqref{eq:constant} is posed as an open problem in Maz'ya's collection \cite[Section~3.11]{Mazya75}. The case $a=0$ is the critical line. The condition at the origin in $\cA_{p,\gamma}$ is already implied by finite energy. Indeed, if $a\le0$ and $\nabla\psi(0)\ne0$, then the gradient integral has nonintegrable local behaviour $\int_0^\varepsilon r^{pa-1}\dd r$. More singular weights may force further Taylor coefficients to vanish. In \cref{lem:annular} we show that the infimum is unchanged when the tests are required to vanish near the origin. This allows us to work on the Emden--Fowler cylinder without altering the original problem.

The case $p=2$ provides a useful starting point. A quarter-turn takes $\nabla\psi$ to a divergence-free field and preserves both norms. Consequently, the planar theorem of Costin and Maz'ya \cite[Theorem~2]{CostinMazya} gives
\begin{equation}\label{eq:hilbert}
 C(2,\gamma)=
 \begin{cases}
 \displaystyle\frac{1+(1-\gamma)^2}{\gamma^2\bigl(3+(1-\gamma)^2\bigr)},
 &-1-\sqrt3<\gamma<-1+\sqrt3,\quad\gamma\ne0,\\[6pt]
 \displaystyle\frac1{1+\gamma^2},
 &\gamma\le-1-\sqrt3\ \text{or}\ \gamma\ge-1+\sqrt3.
 \end{cases}
\end{equation}
The expressions agree at the two endpoints. We recall the correspondence of the admissible classes in \cref{rem:rotation}. Sharp curl-free inequalities and remainder terms were subsequently obtained by Hamamoto and Takahashi \cite{HTsharp,HTremainder}, with a weighted Rellich--Hardy extension in \cite{HTrellich}. In dimensions $N\ge3$, Hamamoto \cite{HamamotoSolenoidal} determined the sharp Hardy--Leray constant for solenoidal fields without symmetry assumptions. Related second-order Caffarelli--Kohn--Nirenberg inequalities were studied in \cite{CFL,DN2025}. Sharp second-order uncertainty principles were obtained by Cazacu, Flynn and Lam \cite{CFLuncertainty}. The recent paper \cite{DN2026} develops weighted $L^2$ inequalities and stability estimates for curl-free fields and second derivatives.

For $p\ne2$, the nonlinear powers couple the angular frequencies, so the sharp quotient no longer separates into Fourier modes. The aim of this paper is to determine its behaviour near the critical line and to identify regions where the radial value remains sharp for arbitrary tests. These are two different parts of the same problem. Near $a=0$, slowly varying affine profiles give the leading term. For large positive or negative $a$, the radial branch is exactly sharp, not only asymptotically sharp.

Our first main result gives the critical coefficient for the full quotient, including $p=1$.
\begin{theorem}\label{thm:critical}
For every fixed $1\le p<\infty$, as $a\to0$,
\begin{equation}\label{eq:critical}
 \Lambda_p(a)=2^{p/2}|a|^p+O_p(|a|^{p+1}),
 \qquad
 |a|^pC(p,\gamma)=2^{-p/2}+O_p(|a|).
\end{equation}
Moreover, $\Lambda_p(0)=0$, or equivalently $C(p,1-2/p)=\infty$.
\end{theorem}

The coefficient in \eqref{eq:critical} is already contained in \eqref{eq:hilbert} when $p=2$. At $p=1$, the theorem gives $\Lambda_1(a)=\sqrt2\,|a|+O(|a|^2)$. The failure at the critical line concerns the whole-plane inequality without a logarithmic correction; the logarithmically corrected $p=2$ inequality on a disk was studied in \cite{HTcritical}.

The reason for the first angular harmonics is natural. Under the transform
\[
 t=\log r,\qquad \psi(r,\theta)=r^\beta\phi(t,\theta),
\]
a vector-valued function $u:\R\to\R^2$ gives the first-harmonic profile $\phi(t,\theta)=u(t)\cdot(\cos\theta,\sin\theta)$, which represents
\[
 \psi(x)=|x|^{-a}u(\log|x|)\cdot x.
\]
At $a=0$, a constant $u$ gives an affine function, whose Hessian vanishes. The proof separates this sector from the other angular modes. A direct angular Poincar\'e estimate controls the complement, while a dual test gives the sharp coefficient on the first harmonics. The same estimates show that every normalized family with transformed Hessian energy $O(|a|^p)$ approaches this angular sector at rate $O(|a|)$; the precise supporting statement is \cref{prop:concentration}.

We next consider the radial branch. Let $C_{\rad}(p,\gamma)$ denote the best constant obtained by restricting \eqref{eq:constant} to radial tests. We prove in \cref{prop:radial} that
\[
 C_{\rad}(p,\gamma)=B^{-p/2}
\]
for every $a\in\R$. The question is when this value remains sharp without a radiality assumption. To state our exterior regions, let $A_0=3.3553013976\ldots$ be the unique root $A>1$ of
\begin{equation}\label{eq:A0}
 A^4-2A^3-3A^2-4A-4=0,
\end{equation}
and put
\begin{equation}\label{eq:thresholds}
 P_p(s)=\bigl[ps^2+p(p-1)s+4\bigr]^2-2p^4s^2-32,
 \qquad E_p^+=\sqrt{2p^2+8}-(p+1).
\end{equation}
For $1\le p<2$, let $A_p^+$ be the largest positive zero of $A\mapsto P_p(A)$. For $p>2$, let $A_p^-$ be the largest positive zero of $A\mapsto P_p(-A)$. These roots exist, as shown in Appendix~\ref{app:polynomial}. In particular, $P_1(s)=(s^2-2)(s^2+8)$, so $A_1^+=\sqrt2$.

\begin{theorem}\label{thm:radiality}
The unrestricted constant satisfies $C(p,\gamma)=B^{-p/2}$ in the following regions:
\begin{equation}\label{eq:rays}
 \begin{aligned}
 1\le p<2:&\qquad a\le-A_0\quad\text{or}\quad a\ge A_p^+,\\
 p=2:&\qquad a\le-1-\sqrt3\quad\text{or}\quad a\ge-1+\sqrt3,\\
 p>2:&\qquad a\le-A_p^-\quad\text{or}\quad a\ge E_p^+.
 \end{aligned}
\end{equation}
\end{theorem}

For $p=2$, the two endpoints are exact by \eqref{eq:hilbert}. For $p\ne2$, \cref{thm:radiality} gives sufficient exterior regions, rather than the exact transition points. Its proof has two parts. An entropy identity gives the positive ray for $p>2$. A homogeneous calibration gives the other rays and also applies at $p=1$. Thus, for each fixed $p\ne2$, only a bounded interval of $a$ remains to be determined.

Inside the transition region, angular tests give the following strict comparison.
\begin{theorem}\label{thm:breaking}
Let $1\le p<\infty$ and $a\ne0$.
\begin{enumerate}[label=\textup{(\roman*)}]
 \item If $-1\le a\le1/2$, then $C(p,\gamma)>B^{-p/2}$.
 \item If $2<p\le4$ and $-1-\sqrt3\le a\le-1+\sqrt3$, then the same strict inequality holds. For $p=2$, it holds in the open interval, whereas equality holds at its two endpoints.
\end{enumerate}
\end{theorem}

The distinction between linear instability and global symmetry breaking is familiar from the classical Caffarelli--Kohn--Nirenberg problem; see \cite{CKN,FelliSchneider,DETT}. Here the transformed gradient also satisfies a differential constraint. In \cref{thm:variation}, we give an explicit necessary and sufficient condition for nonnegativity of the periodic second variation at the radial state. The first angular harmonic controls this condition, but its axial frequency need not be zero. A negative second variation gives strict symmetry breaking. The converse fails: at $(p,a)=(4,7/10)$ the periodic second variation is nonnegative, although the unrestricted constant is strictly larger than the radial value; see \cref{cor:nonnecessity}.

The paper is organised as follows. Section~\ref{sec:cylinder} gives the exact cylinder reduction and the radial constant. In Section~\ref{sec:critical}, we prove the critical asymptotic and the angular concentration estimate. Section~\ref{sec:radiality} proves the exterior radiality regions. Section~\ref{sec:breaking} treats symmetry breaking and the periodic second variation. Section~\ref{sec:interaction} gives an example in which two angular frequencies improve every single-frequency profile, and Section~\ref{sec:remaining} describes the remaining sharp-constant problem. The polynomial verifications and an axial-frequency example are collected in the Appendix.

\section{Exact reduction to the cylinder}\label{sec:cylinder}
We first remove the behaviour at the origin from the variational problem. We then pass to two equivalent formulations on the cylinder: compatible vector fields and scalar potentials. The former is useful for radiality, and the latter for the critical asymptotic and periodic perturbations.

\subsection{Reduction to annular tests}
We call a function \emph{flattened near the origin} if it is constant in a neighbourhood of the origin, and \emph{annular} if it vanishes both near the origin and outside a compact set.
\begin{lemma}\label{lem:annular}
The infimum defining $\Lambda_p(a)$ is unchanged if one first restricts to functions flattened near the origin and then to annular functions. Both reductions preserve radiality.
\end{lemma}
\begin{proof}
Fix $\psi\in\cA_{p,\gamma}$ and put $c=\psi(0)$. Suppose first that $\psi-c$ does not vanish to infinite order at the origin, and let $m\ge1$ be the degree of its first nonzero Taylor term. If $\psi-c$ vanishes to infinite order, choose an integer $m\ge2$ so large that $a+m-1>0$. In either case, Taylor's theorem gives, near the origin,
\begin{equation}\label{eq:taylor}
 |D^j(\psi-c)(x)|\le C_m|x|^{m-j},\qquad j=0,1,2.
\end{equation}
When $m=1$, the estimate for $j=2$ is understood on $0<|x|<1$; boundedness of $D^2\psi$ is more than enough.

If a genuine term of degree $m$ is present, finite gradient energy forces
\begin{equation}\label{eq:localcondition}
 a+m-1>0.
\end{equation}
Indeed, the gradient of the leading homogeneous polynomial is nonzero on a cone of positive angular measure. On a smaller cone and for sufficiently small $r$, $|\nabla\psi|$ is bounded below by a positive multiple of $r^{m-1}$. The gradient energy on this cone is therefore bounded below by a positive multiple of
\[
 \int_0^\varepsilon r^{p(m-1)}r^{p(\gamma-1)}r\dd r
 =\int_0^\varepsilon r^{p(a+m-1)-1}\dd r,
\]
which is finite exactly under \eqref{eq:localcondition}. In the flat case, the same condition holds by the choice of $m$.

Choose a smooth function $\eta$ on $[0,\infty)$ with $\eta=0$ on $[0,1]$ and $\eta=1$ on $[2,\infty)$, and set
\[
 \psi_\varepsilon(x)=c+\eta(|x|/\varepsilon)(\psi(x)-c).
\]
This function equals $c$ on $|x|\le\varepsilon$ and equals $\psi$ on $|x|\ge2\varepsilon$. On the transition annulus, Leibniz's rule and \eqref{eq:taylor} give
\[
 |\nabla(\psi_\varepsilon-\psi)|\le C\varepsilon^{m-1},
 \qquad |D^2(\psi_\varepsilon-\psi)|\le C\varepsilon^{m-2}.
\]
The corresponding bounds in the inner disk have $r$ in place of $\varepsilon$. Hence
\[
 \begin{aligned}
 I_{p,\gamma}[\psi_\varepsilon-\psi]
 &\le C\int_0^{2\varepsilon}r^{p(m-1)}r^{p(\gamma-1)}r\dd r
 \le C\varepsilon^{p(a+m-1)},\\
 J_{p,\gamma}[\psi_\varepsilon-\psi]
 &\le C\int_0^{2\varepsilon}r^{p(m-2)}r^{p\gamma}r\dd r
 \le C\varepsilon^{p(a+m-1)}.
 \end{aligned}
\]
Both expressions tend to zero. The reverse triangle inequality in the two weighted $L^p$ spaces now gives convergence of both norms. Since $I_{p,\gamma}[\psi]>0$, the denominator stays positive for small $\varepsilon$, and the quotients converge. This proves the first reduction.

Now let $\psi$ be flattened near the origin. Its gradient is supported in an annulus $r_0\le|x|\le R_0$. Choose $\lambda>R_0/r_0$ and put
\[
 \Psi_\lambda(x)=\psi(x)-\psi(\lambda x).
\]
Then $\Psi_\lambda$ is annular, and the derivative supports of its two summands are disjoint. A change of variables gives
\[
 I_{p,\gamma}[\psi(\lambda\,\cdot)]=\lambda^{p(1-a)}I_{p,\gamma}[\psi],
 \qquad
 J_{p,\gamma}[\psi(\lambda\,\cdot)]=\lambda^{p(1-a)}J_{p,\gamma}[\psi].
\]
The two energies of $\Psi_\lambda$ are therefore the corresponding sums, and
\[
 \frac{J_{p,\gamma}[\Psi_\lambda]}{I_{p,\gamma}[\Psi_\lambda]}
 =\frac{(1+\lambda^{p(1-a)})J_{p,\gamma}[\psi]}
 {(1+\lambda^{p(1-a)})I_{p,\gamma}[\psi]}
 =\frac{J_{p,\gamma}[\psi]}{I_{p,\gamma}[\psi]}.
\]
Both constructions preserve radiality, so the last assertion follows as well.
\end{proof}

\subsection{Compatible fields and scalar potentials}
Write $x=r(\cos\theta,\sin\theta)$ and $t=\log r$, and let
\[
 e_r=(\cos\theta,\sin\theta),\qquad e_\theta=(-\sin\theta,\cos\theta).
\]
For a test flattened near the origin, define $v,w$ by
\begin{equation}\label{eq:fieldtransform}
 \nabla\psi(r,\theta)=e^{-at}\bigl(v(t,\theta)e_r+w(t,\theta)e_\theta\bigr).
\end{equation}
Equivalently, $\dd\psi=e^{\beta t}(v\dd t+w\dd\theta)$. This one-form is closed, so
\begin{equation}\label{eq:compatible}
 w_t-v_\theta+\beta w=0.
\end{equation}
We call a smooth pair satisfying \eqref{eq:compatible} \emph{compatible}.

Since $\partial_r=e^{-t}\partial_t$, differentiation gives
\[
 \begin{aligned}
 \psi_{rr}&=e^{-(a+1)t}(v_t-av),\\
 r^{-1}\psi_{r\theta}-r^{-2}\psi_\theta&=e^{-(a+1)t}(v_\theta-w),\\
 r^{-2}\psi_{\theta\theta}+r^{-1}\psi_r&=e^{-(a+1)t}(w_\theta+v).
 \end{aligned}
\]
In the orthonormal polar frame,
\[
 |D^2\psi|_{\Frob}^2
 =\psi_{rr}^2+2(r^{-1}\psi_{r\theta}-r^{-2}\psi_\theta)^2
 +(r^{-2}\psi_{\theta\theta}+r^{-1}\psi_r)^2.
\]
The powers of $r$ cancel against the weights and the Jacobian. Thus the two planar energies become
\begin{align*}
 I_a[v,w]&=\int_{\R\times S^1}(v^2+w^2)^{p/2}\dd t\dd\theta,
 \\
 J_a[v,w]&=\int_{\R\times S^1}
 \bigl((v_t-av)^2+2(v_\theta-w)^2+(w_\theta+v)^2\bigr)^{p/2}\dd t\dd\theta.
\end{align*}
\begin{proposition}\label{prop:fields}
For every $1\le p<\infty$ and $a\in\R$, $\Lambda_p(a)$ is the infimum of $J_a[v,w]/I_a[v,w]$ over all nonzero smooth compactly supported compatible pairs.
\end{proposition}
\begin{proof}
The formulas above send every flattened planar test to a compactly supported compatible pair without changing its quotient. Conversely, let $(v,w)$ be such a pair. The one-form
\[
 \omega=e^{\beta t}(v\dd t+w\dd\theta)
\]
is closed by \eqref{eq:compatible}. The only possible obstruction to exactness on $\R\times S^1$ is its period around an angular circle. Denote this period by
\[
 p(t)=\int_0^{2\pi}e^{\beta t}w(t,\theta)\dd\theta.
\]
Compatibility and periodicity imply
\[
 p'(t)=\int_0^{2\pi}e^{\beta t}(w_t+\beta w)\dd\theta
 =\int_0^{2\pi}e^{\beta t}v_\theta\dd\theta=0.
\]
The period vanishes for large $|t|$, because the pair is compactly supported. Hence it vanishes identically, and $\omega$ has a global smooth potential.

This potential is constant near both ends of the cylinder, although the constants may differ. Subtract its value near $+\infty$ and then replace $t$ by $\log r$. The resulting planar function vanishes for large $r$ and is constant for small $r$. It extends smoothly across the origin, and \eqref{eq:fieldtransform} recovers the given pair. Its gradient has compact annular support, so both weighted energies are finite for every $a$, including $a=0$. The conclusion now follows from \cref{lem:annular}.
\end{proof}

\begin{remark}\label{rem:rotation}
Let $R$ denote rotation by $\pi/2$. If $\psi\in C_c^\infty(\R^2)$, then $R\nabla\psi$ is divergence free and
\[
 |R\nabla\psi|=|\nabla\psi|,\qquad
 |\nabla(R\nabla\psi)|_{\Frob}=|D^2\psi|_{\Frob}.
\]
Conversely, if $U\in C_c^\infty(\R^2;\R^2)$ is divergence free, then $-RU$ is curl free and hence equals $\nabla\psi$ for a smooth potential on $\R^2$. The potential is constant outside a large ball; subtracting that constant makes it compactly supported. Moreover, $U(0)=0$ if and only if $\nabla\psi(0)=0$. Together with \cref{lem:annular}, this identifies the admissible classes in the planar divergence-free inequality and in \eqref{eq:constant} for $p=2$.
\end{remark}

For the scalar formulation, take an annular test and write
\begin{equation}\label{eq:scalartransform}
 \psi(r,\theta)=r^\beta\phi(\log r,\theta).
\end{equation}
Then $\psi_r=r^{\beta-1}(\phi_t+\beta\phi)$ and $r^{-1}\psi_\theta=r^{\beta-1}\phi_\theta$. Accordingly, define
\begin{align}
 P_a\phi&=(\phi_t+\beta\phi,\phi_\theta),\label{eq:Pa}\\
 Q_a\phi&=\bigl((\partial_t-a)(\partial_t+\beta)\phi,
 \sqrt2(\partial_t-a)\phi_\theta,
 \phi_{\theta\theta}+\phi_t+\beta\phi\bigr).
 \label{eq:Qa}
\end{align}
Every $\phi\in C_c^\infty(\R\times S^1)$ gives an annular planar test through \eqref{eq:scalartransform}, and every annular test has such a representative. Therefore
\begin{equation}\label{eq:scalarinf}
 \Lambda_p(a)=\inf_{0\ne\phi\in C_c^\infty(\R\times S^1)}
 \frac{\displaystyle\int_{\R\times S^1}|Q_a\phi|^p}
 {\displaystyle\int_{\R\times S^1}|P_a\phi|^p}.
\end{equation}
The denominator is positive for nonzero $\phi$: the equations $P_a\phi=0$ imply $\phi_\theta=0$ and $\phi_t+\beta\phi=0$, and compact support then gives $\phi=0$. The equality in \eqref{eq:scalarinf} concerns the infima. A compactly supported compatible pair can have a planar potential with different constants at the two cylinder ends, so the pair need not equal $P_a\phi$ for a compactly supported $\phi$. The annular reduction is what identifies the two infima.

We shall also use a fixed Cartesian frame. Put
\[
 G(t,\theta)=v(t,\theta)e_r+w(t,\theta)e_\theta.
\]
The field is called compatible when its polar components satisfy \eqref{eq:compatible}. Differentiating the moving frame gives
\[
 G_t-aG=(v_t-av)e_r+(w_t-aw)e_\theta,
 \qquad
 G_\theta=(v_\theta-w)e_r+(w_\theta+v)e_\theta.
\]
Since $\beta+a=1$, compatibility yields $w_t-aw=v_\theta-w$. Consequently,
\begin{equation}\label{eq:cartesian}
 I_a=\int|G|^p,
 \qquad
 J_a=\int\bigl(|G_t-aG|^2+|G_\theta|^2\bigr)^{p/2}.
\end{equation}
Unless another domain is indicated, cylinder integrals and norms use the measure $\dd t\dd\theta$ on $\R\times S^1$.

\subsection{Periodic and long cutoffs}
The next two lemmas allow us to use periodic or purely angular profiles in \eqref{eq:scalarinf}. For $L>0$, let $\T_L=\R/L\mathbb Z$.
\begin{lemma}\label{lem:periodic}
Let $\phi\in C^\infty(\T_L\times S^1)$ and suppose $\int_{\T_L\times S^1}|P_a\phi|^p>0$. There are $\phi_N\in C_c^\infty(\R\times S^1)$ such that
\[
 \frac{\int|Q_a\phi_N|^p}{\int|P_a\phi_N|^p}
 \longrightarrow
 \frac{\int_{\T_L\times S^1}|Q_a\phi|^p}
 {\int_{\T_L\times S^1}|P_a\phi|^p}.
\]
\end{lemma}
\begin{proof}
Extend $\phi$ periodically. Choose $\zeta_N\in C_c^\infty(\R)$ equal to one on $[0,NL]$, supported in $[-1,NL+1]$, with its first two derivatives bounded independently of $N$. Put $\phi_N=\zeta_N\phi$. Product differentiation gives
\[
 P_a(\zeta_N\phi)=\zeta_NP_a\phi+(\zeta_N'\phi,0)
\]
and
\[
 \begin{aligned}
 (Q_a(\zeta_N\phi))_1&=\zeta_N(Q_a\phi)_1+2\zeta_N'\phi_t
 +\zeta_N''\phi+(\beta-a)\zeta_N'\phi,\\
 (Q_a(\zeta_N\phi))_2&=\zeta_N(Q_a\phi)_2+\sqrt2\,\zeta_N'\phi_\theta,\\
 (Q_a(\zeta_N\phi))_3&=\zeta_N(Q_a\phi)_3+\zeta_N'\phi.
 \end{aligned}
\]
The $N$ complete cells contribute exactly $N$ times the two periodic energies. The remaining terms are supported in two strips of fixed width. The profile, its derivatives, and the cutoff derivatives are uniformly bounded there. Thus
\[
 \int|P_a\phi_N|^p=N\int_{\T_L\times S^1}|P_a\phi|^p+O(1),
 \qquad
 \int|Q_a\phi_N|^p=N\int_{\T_L\times S^1}|Q_a\phi|^p+O(1).
\]
The leading denominator is positive. Dividing proves the assertion.
\end{proof}

For $\Phi\in C^\infty(S^1)$, define
\begin{equation}\label{eq:angularoperators}
 P_a^{\ang}\Phi=(\beta\Phi,\Phi'),
 \qquad
 Q_a^{\ang}\Phi=(-a\beta\Phi,-\sqrt2a\Phi',\Phi''+\beta\Phi).
\end{equation}
\begin{lemma}\label{lem:long}
Suppose $P_a^{\ang}\Phi\not\equiv0$. If $0\ne\chi\in C_c^\infty(\R)$ and $\phi_R(t,\theta)=\chi(t/R)\Phi(\theta)$, then
\[
 \begin{aligned}
 \frac1R\int|P_a\phi_R|^p&\longrightarrow
 \left(\int_\R|\chi|^p\right)\left(\int_0^{2\pi}|P_a^{\ang}\Phi|^p\dd\theta\right),\\
 \frac1R\int|Q_a\phi_R|^p&\longrightarrow
 \left(\int_\R|\chi|^p\right)\left(\int_0^{2\pi}|Q_a^{\ang}\Phi|^p\dd\theta\right).
 \end{aligned}
\]
In particular,
\begin{equation}\label{eq:angularbound}
 \Lambda_p(a)\le
 \frac{\int_0^{2\pi}|Q_a^{\ang}\Phi|^p\dd\theta}
 {\int_0^{2\pi}|P_a^{\ang}\Phi|^p\dd\theta}.
\end{equation}
\end{lemma}
\begin{proof}
After $t=Rs$, the fields are
\[
 P_a\phi_R=\chi(s)P_a^{\ang}\Phi+R^{-1}(\chi'(s)\Phi,0)
\]
and
\[
 \begin{aligned}
 Q_a\phi_R={}&\chi(s)Q_a^{\ang}\Phi
 +R^{-1}\bigl((\beta-a)\chi'(s)\Phi,\sqrt2\chi'(s)\Phi',\chi'(s)\Phi\bigr)\\
 &+R^{-2}(\chi''(s)\Phi,0,0).
 \end{aligned}
\]
All terms are supported on the fixed set $\operatorname{supp}\chi\times S^1$ and are uniformly bounded for $R\ge1$. Dominated convergence proves the two limits. The first limit is positive, so division and \eqref{eq:scalarinf} give \eqref{eq:angularbound}.
\end{proof}

\subsection{The radial branch and universal bounds}
We first note a zero-set form of the chain rule, which will be useful to justify the endpoint arguments later.
\begin{lemma}\label{lem:zeros}
Let $\Omega\subset\R^d$ be open and $U:\Omega\to\R^m$ be locally Lipschitz. Then $DU=0$ almost everywhere on $\{U=0\}$. If $F:\R^m\to\R$ is locally Lipschitz and continuously differentiable away from the origin, then
\[
 D(F\circ U)=DF(U)DU
\]
almost everywhere on $\{U\ne0\}$. With the right-hand side set equal to zero on $\{U=0\}$, the identity holds almost everywhere in $\Omega$. Here, ``almost everywhere'' means with respect to the ambient Lebesgue measure on $\Omega$.
\end{lemma}
\begin{proof}
For a scalar locally Lipschitz function, the restriction to almost every line parallel to a coordinate axis is absolutely continuous. Its derivative vanishes almost everywhere on each level set. Fubini's theorem therefore gives $Df=0$ almost everywhere on $\{f=c\}$. Applying this fact to the components of $U$ proves the first assertion. On $\{U\ne0\}$ the ordinary chain rule applies. On $\{U=0\}$, both $DU$ and $D(F\circ U)$ vanish almost everywhere, by the same level-set argument. The composition is locally Lipschitz, since both maps are locally Lipschitz. The result also holds on the cylinder by using coordinate charts.
\end{proof}

A radial compatible field has the form $(v,w)=(h,0)$, so its energies are
\begin{equation*}
 I_{\rad}[h]=2\pi\int_\R|h|^p\dd t,
 \qquad
 J_{\rad}[h]=2\pi\int_\R\bigl((h'-ah)^2+h^2\bigr)^{p/2}\dd t.
\end{equation*}
\begin{proposition}\label{prop:radial}
For every $h\in C_c^\infty(\R)$,
\[
 J_{\rad}[h]\ge B^{p/2}I_{\rad}[h],
\]
and the constant is sharp. Consequently, $C_{\rad}(p,\gamma)=B^{-p/2}$ for every $a\in\R$.
\end{proposition}
\begin{proof}
On $\{h\ne0\}$, project $(h'-ah,h)$ onto the unit vector
\[
 \sgn(h)\left(-\frac a{\sqrt B},\frac1{\sqrt B}\right).
\]
Set this vector equal to zero on $\{h=0\}$. By \cref{lem:zeros}, $h'=0$ almost everywhere on that set, and hence
\[
 \bigl((h'-ah)^2+h^2\bigr)^{1/2}
 \ge-\frac a{\sqrt B}\sgn(h)h'+\sqrt B\,|h|
\]
almost everywhere. Multiply by $|h|^{p-1}$ and integrate. The derivative term vanishes because
\[
 \sgn(h)h'|h|^{p-1}=\frac1p(|h|^p)'
\]
almost everywhere, also at $p=1$. If $p>1$, H\"older's inequality gives
\[
 \sqrt B\int_\R|h|^p
 \le\left(\int_\R\bigl((h'-ah)^2+h^2\bigr)^{p/2}\right)^{1/p}
 \left(\int_\R|h|^p\right)^{(p-1)/p}.
\]
For nonzero $h$, division and taking the $p$th power prove the claimed estimate. For $p=1$, the integrated projection inequality is already the desired estimate. The zero function causes no exception.

For sharpness, choose $0\ne\eta\in C_c^\infty(\R)$ and let $h_R(t)=\eta(t/R)$. Changing variables and applying dominated convergence gives
\[
 \frac{J_{\rad}[h_R]}{I_{\rad}[h_R]}\longrightarrow B^{p/2}.
\]
These fields come from admissible radial planar functions. Indeed, put
\[
 f_R(r)=r^{-a}h_R(\log r),\qquad
 \Psi_R(r)=-\int_r^\infty f_R(s)\dd s.
\]
The function $f_R$ has compact annular support, so $\Psi_R$ is smooth, constant near the origin, and zero for large $r$. Its transformed radial field is $(h_R,0)$. The radial reduction in \cref{lem:annular} now gives the stated value of $C_{\rad}$.
\end{proof}

The following shifted-derivative estimate will be used both for the universal bound and for the first angular harmonics.
\begin{proposition}\label{prop:shift}
Let $u\in C_c^\infty(\R;\R^m)$ and $a\in\R$. Then
\begin{equation}\label{eq:shift}
 \int_\R|u'-au|^p\dd t\ge|a|^p\int_\R|u|^p\dd t.
\end{equation}
Consequently,
\begin{equation}\label{eq:universal}
 |a|^p\le\Lambda_p(a)\le B^{p/2}.
\end{equation}
For $a\ne0$, this gives $B^{-p/2}\le C(p,\gamma)\le|a|^{-p}$ and $|a|^pC(p,\gamma)\to1$ as $|a|\to\infty$.
\end{proposition}
\begin{proof}
Suppose first that $p>1$ and $u\not\equiv0$. The chain rule and compact support give
\[
 \int_\R|u|^{p-2}u\cdot u'=\frac1p\int_\R(|u|^p)'=0,
\]
where $|u|^{p-2}u$ is zero at $u=0$. Therefore
\[
 |a|\int_\R|u|^p
 =\left|\int_\R|u|^{p-2}u\cdot(u'-au)\right|
 \le\left(\int_\R|u|^p\right)^{(p-1)/p}
 \left(\int_\R|u'-au|^p\right)^{1/p}.
\]
Division proves \eqref{eq:shift}. If $p=1$, define $u/|u|=0$ on its zero set. By \cref{lem:zeros},
\[
 \int_\R\frac u{|u|}\cdot u'=\int_\R(|u|)'=0.
\]
It follows directly that
\[
 |a|\int_\R|u|=\left|\int_\R\frac u{|u|}\cdot(u'-au)\right|
 \le\int_\R|u'-au|.
\]
The case $u=0$ is immediate.

Apply \eqref{eq:shift} to $G(\,\cdot\,,\theta)$ for each fixed $\theta$ and integrate in $\theta$. Then
\[
 |a|^p\int|G|^p\le\int|G_t-aG|^p
 \le\int\bigl(|G_t-aG|^2+|G_\theta|^2\bigr)^{p/2}.
\]
For compatible fields, \eqref{eq:cartesian} and \cref{prop:fields} give the lower bound in \eqref{eq:universal}. The radial sequence gives the upper bound. Taking reciprocals gives the bounds for $C$. Finally,
\[
 \frac{|a|^p}{B^{p/2}}\le |a|^pC(p,\gamma)\le1,
\]
and the left-hand side tends to one as $|a|\to\infty$.
\end{proof}

\section{The critical asymptotic and angular concentration}\label{sec:critical}
We now prove \cref{thm:critical}. The first angular harmonics play a special role because they represent affine functions at the critical parameter. We first control the other harmonics by an angular Poincar\'e inequality. We then use a dual test on the first harmonics to obtain the sharp leading coefficient.

\subsection{The first angular harmonics and their complement}
Write
\[
 e(\theta)=(\cos\theta,\sin\theta),\qquad
 e'(\theta)=(-\sin\theta,\cos\theta),
\]
and define
\begin{equation}\label{eq:projection}
 \Pi\phi(t,\theta)=u(t)\cdot e(\theta),
 \qquad u(t)=\frac1\pi\int_0^{2\pi}\phi(t,\theta)e(\theta)\dd\theta.
\end{equation}
For a vector-valued $u$ on $\R$, put
\[
 Wu(t,\theta)=\bigl(u(t)\cdot e(\theta),u(t)\cdot e'(\theta)\bigr).
\]
In particular, $|Wu(t,\theta)|=|u(t)|$.

The projection acts componentwise on vector-valued functions. Its kernel is
\[
 \Pi F(t,\theta)=\frac1\pi\int_0^{2\pi}\cos(\theta-\vartheta)F(t,\vartheta)\dd\vartheta,
\]
so Young's inequality gives
\begin{equation}\label{eq:projectionbound}
 \norm{\Pi F}_{L^p}\le\frac4\pi\norm{F}_{L^p},\qquad 1\le p\le\infty.
\end{equation}
Both $P_a$ and $Q_a$ commute with $\Pi$.

\begin{lemma}\label{lem:complement}
For each $1\le p<\infty$, there is $C_p<\infty$ such that
\begin{equation*}
 \norm{P_ar}_{L^p}\le C_p\norm{Q_ar}_{L^p}
\end{equation*}
for every $a\in\R$ and every $r\in C_c^\infty(\R\times S^1)$ satisfying $\Pi r=0$. The constant is independent of $a$.
\end{lemma}
\begin{proof}
Consider the associated Cartesian field
\[
 H=(r_t+\beta r)e_r+r_\theta e_\theta.
\]
Since $\Pi r=0$, one has $\int_0^{2\pi}r e_r\dd\theta=0$ for each $t$, and the same is true with $r_t$ in place of $r$. Integration by parts also gives
\[
 \int_0^{2\pi}r_\theta e_\theta\dd\theta
 =-\int_0^{2\pi}r\,\partial_\theta e_\theta\dd\theta
 =\int_0^{2\pi}r e_r\dd\theta=0.
\]
Thus $H(t,\,\cdot\,)$ has zero angular mean. Differentiating the moving frame, and using $\beta-1=-a$, we obtain
\[
 H_\theta=(\partial_t-a)r_\theta\,e_r
 +(r_{\theta\theta}+r_t+\beta r)e_\theta.
\]
The second and third components of \eqref{eq:Qa} therefore give
\[
 |H_\theta|^2=\frac12|(Q_ar)_2|^2+|(Q_ar)_3|^2\le|Q_ar|^2.
\]

For completeness, a zero-mean periodic vector function $h$ satisfies
\[
 h(\theta)=\frac1{2\pi}\int_0^{2\pi}\bigl(h(\theta)-h(\vartheta)\bigr)\dd\vartheta.
\]
The fundamental theorem of calculus bounds each difference by $\int_0^{2\pi}|h'|$. Hence
\[
 \norm{h}_{L^p(S^1)}\le2\pi\norm{h'}_{L^p(S^1)},\qquad1\le p<\infty,
\]
where H\"older's inequality is used when $p>1$, and the same estimate follows directly when $p=1$. Applying this to $H(t,\,\cdot\,)$ and integrating in $t$ proves
\[
 \norm{P_ar}_{L^p}=\norm{H}_{L^p}
 \le2\pi\norm{H_\theta}_{L^p}\le2\pi\norm{Q_ar}_{L^p}.
\]
No restriction on $a$ was used.
\end{proof}

The next proposition records the resulting angular localization. The hypothesis does not require the family to be minimizing.
\begin{proposition}\label{prop:concentration}
Fix $1\le p<\infty$. For $\phi\in C_c^\infty(\R\times S^1)$, write $\Pi\phi=u\cdot e$ and $r=(I-\Pi)\phi$. There is a constant $C_p$, independent of $a\in\R$ and $\phi$, such that
\begin{equation}\label{eq:concentration}
 \norm{P_ar}_{L^p}+\norm{u'-au}_{L^p(\R)}
 +\norm{P_a\phi-Wu}_{L^p}\le C_p\norm{Q_a\phi}_{L^p}.
\end{equation}
Consequently, if $M>0$ is fixed and
\[
 \norm{P_a\phi}_{L^p}=1,\qquad \norm{Q_a\phi}_{L^p}\le M|a|,
\]
then each term on the left of \eqref{eq:concentration} is $O_{p,M}(|a|)$, and, for sufficiently small $|a|$,
\begin{equation}\label{eq:mass}
 \left|2\pi\int_\R|u|^p\dd t-1\right|\le C_{p,M}|a|.
\end{equation}
\end{proposition}
\begin{proof}
By commutation and \eqref{eq:projectionbound},
\[
 \norm{Q_ar}_{L^p}\le\left(1+\frac4\pi\right)\norm{Q_a\phi}_{L^p}.
\]
Combining this with \cref{lem:complement} bounds the first term in \eqref{eq:concentration}.

Put $b=u'-au$. Since $\partial_{\theta\theta}(u\cdot e)=-u\cdot e$, the third component of $Q_a\Pi\phi$ is $b\cdot e$. Rotational invariance gives
\[
 \norm{b\cdot e}_{L^p(\R\times S^1)}
 =\left(\int_0^{2\pi}|\cos\theta|^p\dd\theta\right)^{1/p}\norm{b}_{L^p(\R)}.
\]
The angular factor is positive. Also,
\[
 \norm{b\cdot e}_{L^p}\le\norm{Q_a\Pi\phi}_{L^p}
 \le\frac4\pi\norm{Q_a\phi}_{L^p}.
\]
This proves the bound for $b$. Finally, direct substitution into \eqref{eq:Pa} gives
\[
 P_a\Pi\phi=\bigl((u+b)\cdot e,u\cdot e'\bigr),
 \qquad P_a\Pi\phi-Wu=(b\cdot e,0).
\]
Together with $P_a\phi=P_ar+P_a\Pi\phi$, this proves \eqref{eq:concentration}.

Under the normalized low-energy hypothesis, the reverse triangle inequality yields
\[
 \bigl|\norm{Wu}_{L^p}-1\bigr|\le C_{p,M}|a|.
\]
For small enough $|a|$, we have $1/2\le\norm{Wu}_{L^p}\le3/2$. The map $s\mapsto s^p$ is Lipschitz on this interval. Since
\[
 \norm{Wu}_{L^p}^p=2\pi\int_\R|u|^p\dd t,
\]
we obtain \eqref{eq:mass}.
\end{proof}

The proposition identifies the angular sector of a low-energy family, but it leaves its axial profile $u(t)$ free. It is therefore an angular localization estimate rather than a classification of near-extremizers. For the sharp critical law, we still need the precise estimate of the shifted derivative $u'-au$.

\subsection{A lower bound with the sharp critical coefficient}
For $\Pi\phi=u\cdot e$, put
\[
 b=u'-au,\qquad c=b'+\beta b.
\]
The transformed Hessian on this sector is
\begin{equation}\label{eq:resonantQ}
 Q_a\Pi\phi=\bigl(c\cdot e,\sqrt2 b\cdot e',b\cdot e\bigr).
\end{equation}
Although the nonlinear energy does not split into angular modes, a suitable dual field can eliminate the complementary modes from the pairing.

\begin{lemma}\label{lem:duality}
Let $1\le p<\infty$ and $|a|<1$. Then every $\phi\in C_c^\infty(\R\times S^1)$ satisfies
\begin{equation}\label{eq:dualbound}
 \int|Q_a\phi|^p\ge2\pi\,2^{p/2}\left(1-\frac a4\right)^p
 \int_\R|u'-au|^p\dd t,
\end{equation}
where $u$ is given by \eqref{eq:projection}.
\end{lemma}
\begin{proof}
Define
\[
 Z_b=(b\cdot e,\sqrt2 b\cdot e',b\cdot e),
 \qquad Y=|Z_b|^{p-2}Z_b,
\]
with $Y=0$ on $\{b=0\}$. We observe that
\[
 |Z_b|^2=2(b\cdot e)^2+2(b\cdot e')^2=2|b|^2.
\]
Thus its length is independent of $\theta$. Multiplication by $|Z_b|^{p-2}$ does not introduce other angular frequencies: every component of $Y$ remains in the first-harmonic sector. Since $Q_a$ commutes with $\Pi$, it follows that
\[
 \int Q_a(I-\Pi)\phi\cdot Y=0.
\]
All expressions containing $|b|^{p-2}b$ below are defined to be zero on $\{b=0\}$.

Using \eqref{eq:resonantQ} and the identities
\[
 \int_0^{2\pi}(r\cdot e)(s\cdot e)\dd\theta=\pi r\cdot s,
 \qquad
 \int_0^{2\pi}(r\cdot e')^2\dd\theta=\pi|r|^2,
\]
we find that
\[
 \begin{aligned}
 \int_0^{2\pi}Q_a\Pi\phi\cdot Y\dd\theta
 &=\pi2^{(p-2)/2}\bigl(|b|^{p-2}b\cdot c+3|b|^p\bigr)\\
 &=\pi2^{(p-2)/2}\bigl(|b|^{p-2}b\cdot b'+(4-a)|b|^p\bigr).
 \end{aligned}
\]
The first term integrates to zero in $t$, because
\[
 |b|^{p-2}b\cdot b'=\frac1p(|b|^p)'
\]
almost everywhere. For $p=1$ this is the norm chain rule in \cref{lem:zeros}; compact support removes the boundary term. We have therefore proved
\begin{equation}\label{eq:dualpair}
 \int Q_a\phi\cdot Y
 =\pi2^{(p-2)/2}(4-a)\int_\R|b|^p\dd t.
\end{equation}

If $b=0$ almost everywhere, the conclusion is immediate. Otherwise, for $p>1$ put $p'=p/(p-1)$. Then
\[
 \int|Y|^{p'}=\int|Z_b|^p
 =2\pi\,2^{p/2}\int_\R|b|^p\dd t.
\]
H\"older's inequality in \eqref{eq:dualpair}, followed by division and taking the $p$th power, gives
\[
 \int|Q_a\phi|^p
 \ge2\pi\,2^{p/2}\left(\frac{4-a}{4}\right)^p\int_\R|b|^p\dd t.
\]
For $p=1$, $|Y|\le1$, so \eqref{eq:dualpair} directly gives
\[
 \norm{Q_a\phi}_{L^1}
 \ge\frac{\pi(4-a)}{\sqrt2}\int_\R|b|\dd t
 =2\pi\sqrt2\left(1-\frac a4\right)\int_\R|b|\dd t.
\]
These are exactly \eqref{eq:dualbound}.
\end{proof}

\subsection{Proof of the critical asymptotic}
\begin{proof}[Proof of \cref{thm:critical}]
Take $\Phi(\theta)=\cos\theta$ in \cref{lem:long}. From \eqref{eq:angularoperators},
\[
 \Lambda_p(a)\le |a|^p
 \frac{\displaystyle\int_0^{2\pi}
 \bigl(((1-a)^2+1)\cos^2\theta+2\sin^2\theta\bigr)^{p/2}\dd\theta}
 {\displaystyle\int_0^{2\pi}
 \bigl((1-a)^2\cos^2\theta+\sin^2\theta\bigr)^{p/2}\dd\theta}.
\]
For $|a|\le1/2$, both radicands are uniformly positive. The quotient of angular integrals is therefore differentiable with bounded derivative, and its value at $a=0$ is $2^{p/2}$. Hence, for small $a$,
\begin{equation}\label{eq:criticalupper}
 \Lambda_p(a)\le2^{p/2}|a|^p+C_p|a|^{p+1}.
\end{equation}

To prove the reverse estimate, let $a\ne0$ and choose a test whose quotient is within $|a|^{p+1}$ of the infimum. By homogeneity, we can normalize it as
\begin{equation}\label{eq:nearmin}
 \int|P_a\phi_a|^p=1,
 \qquad
 \int|Q_a\phi_a|^p\le\Lambda_p(a)+|a|^{p+1}.
\end{equation}
The upper bound \eqref{eq:criticalupper} gives $\norm{Q_a\phi_a}_{L^p}\le M_p|a|$ for sufficiently small nonzero $a$, with $M_p$ depending only on $p$. Write $\Pi\phi_a=u_a\cdot e$ and $b_a=u_a'-au_a$. By \cref{prop:concentration},
\begin{equation}\label{eq:nearmass}
 2\pi\int_\R|u_a|^p\dd t=1+O_p(|a|).
\end{equation}
Now apply \cref{lem:duality} and then the shifted-derivative estimate \eqref{eq:shift}. We obtain
\[
 \begin{aligned}
 \int|Q_a\phi_a|^p
 &\ge2\pi\,2^{p/2}\left(1-\frac a4\right)^p\int_\R|b_a|^p\\
 &\ge2^{p/2}|a|^p\left(1-\frac a4\right)^p
 \left(2\pi\int_\R|u_a|^p\right)\\
 &\ge2^{p/2}|a|^p-C_p|a|^{p+1}.
 \end{aligned}
\]
The last step uses \eqref{eq:nearmass} and $(1-a/4)^p=1+O_p(|a|)$. By \eqref{eq:nearmin},
\[
 \Lambda_p(a)\ge\int|Q_a\phi_a|^p-|a|^{p+1}
 \ge2^{p/2}|a|^p-C_p|a|^{p+1}.
\]
Together with \eqref{eq:criticalupper}, this proves the first formula in \eqref{eq:critical}. The leading coefficient is positive, so taking reciprocals gives the second formula.

Finally, take $a=0$, choose $0\ne\chi\in C_c^\infty(\R)$, and put
\[
 \phi_R(t,\theta)=\chi(t/R)\cos\theta.
\]
Its planar representative is
\[
 \psi_R(x)=x_1\chi\left(\frac{\log|x|}{R}\right)\quad(x\ne0),
 \qquad\psi_R(0)=0.
\]
It vanishes near the origin and for large $|x|$, so it belongs to the original admissible class. Writing $s=t/R$, we have
\[
 P_0\phi_R=\bigl((\chi(s)+R^{-1}\chi'(s))\cos\theta,-\chi(s)\sin\theta\bigr)
\]
and
\[
 Q_0\phi_R=\bigl((R^{-2}\chi''(s)+R^{-1}\chi'(s))\cos\theta,
 -\sqrt2R^{-1}\chi'(s)\sin\theta,R^{-1}\chi'(s)\cos\theta\bigr).
\]
Dominated convergence after $t=Rs$ gives
\[
 \frac1R\int|P_0\phi_R|^p\longrightarrow2\pi\int_\R|\chi|^p>0.
\]
For $R\ge1$, the Hessian is bounded by $C R^{-1}(|\chi'(s)|+|\chi''(s)|)$, whence
\[
 \int|Q_0\phi_R|^p=O_{p,\chi}(R^{1-p}).
\]
The quotient is $O(R^{-p})$ and tends to zero. Thus $\Lambda_p(0)=0$.
\end{proof}

\section{Sharpness regions for the radial branch}\label{sec:radiality}
To prove \cref{thm:radiality}, it remains to obtain the lower bound $J_a\ge B^{p/2}I_a$ for arbitrary compatible fields in the stated regions. We use two arguments. The first extracts an integral identity from compatibility and applies Jensen's inequality, so it is suited to $p\ge2$. The second constructs a dual field whose norm gives the radial constant, and it also works when $1\le p<2$. The comparison of radial and unrestricted constants is a recurring issue in weighted second-order inequalities; see, for example, \cite{DoLamLu}.

\subsection{The entropy argument}
Let $G=ve_r+we_\theta$ be a smooth compactly supported compatible field and put $R=|G|$. On $U=\{G\ne0\}$, write locally
\[
 G=R(\cos\varphi,\sin\varphi),\qquad
 q_t=\frac{R_t}{R},\qquad q_\theta=\frac{R_\theta}{R},
\]
and let $c=\cos(\varphi-\theta)$ and $s=\sin(\varphi-\theta)$. These quantities do not require a global choice of the phase. Indeed,
\[
 q_j=\frac{G\cdot G_j}{R^2},\qquad
 \varphi_j=\frac{G_1(G_2)_j-G_2(G_1)_j}{R^2},\qquad j=t,\theta,
\]
are globally defined on $U$, and $c=v/R$, $s=w/R$.

Substitute $v=Rc$ and $w=Rs$ into \eqref{eq:compatible}. Since
\[
 w_t=R(sq_t+c\varphi_t),\qquad
 v_\theta=R(cq_\theta-s\varphi_\theta+s),
\]
compatibility becomes
\begin{equation}\label{eq:amplitude}
 cq_\theta-sq_t-c\varphi_t-s\varphi_\theta+as=0.
\end{equation}
We shall use this relation to control the angular energy of $G$.

\begin{lemma}\label{lem:entropy}
For $p\ge2$ and $a\in\R$, every smooth compactly supported compatible field satisfies
\begin{equation}\label{eq:entropy}
 \int R^p\bigl[1+(pa-2)s^2\bigr]
 =\int R^p\bigl[(p-2)sc\,\varphi_t+(c^2+(p-1)s^2)\varphi_\theta\bigr],
\end{equation}
where the integrands are extended by zero on $\{G=0\}$.
\end{lemma}
\begin{proof}
First work on a regular superlevel set $\Omega_\varepsilon=\{R>\varepsilon\}$. Multiply \eqref{eq:amplitude} by $pR^ps$ and integrate. Since $pR^pq_j=\partial_j(R^p)$,
\[
 \begin{aligned}
 p\int_{\Omega_\varepsilon}R^pscq_\theta
 &=-\int_{\Omega_\varepsilon}R^p(sc)_\theta
 +\int_{\partial\Omega_\varepsilon}R^pscn_\theta\dd\mathcal H^1,\\
 -p\int_{\Omega_\varepsilon}R^ps^2q_t
 &=\int_{\Omega_\varepsilon}R^p(s^2)_t
 -\int_{\partial\Omega_\varepsilon}R^ps^2n_t\dd\mathcal H^1.
 \end{aligned}
\]
The identities
\[
 (sc)_\theta=(c^2-s^2)(\varphi_\theta-1),\qquad
 (s^2)_t=2sc\varphi_t
\]
therefore give
\[
 \begin{aligned}
 0={}&\int_{\Omega_\varepsilon}R^p
 \bigl[1+(pa-2)s^2-(p-2)sc\varphi_t
 -(c^2+(p-1)s^2)\varphi_\theta\bigr]\\
 &+\int_{\partial\Omega_\varepsilon}R^p(scn_\theta-s^2n_t)\dd\mathcal H^1.
 \end{aligned}
\]
The absolute value of the boundary term is at most $C\varepsilon^p\operatorname{Per}(\Omega_\varepsilon)$.

We can choose regular levels along which this term tends to zero. Indeed, let $\delta_j\downarrow0$. The coarea formula and Sard's theorem give regular values $\varepsilon_j\in(\delta_j/2,\delta_j)$ for which
\[
 \varepsilon_j^p\operatorname{Per}(\Omega_{\varepsilon_j})
 \le2\delta_j^{p-1}\int_{\{0<R<\delta_j\}}|\nabla R|
 \le2\delta_j^{p-1}\norm{DG}_{L^1}\longrightarrow0.
\]
Here $|\nabla R|\le|DG|$, and $G$ is smooth and compactly supported. The difference between the interior integral over $\Omega_{\varepsilon_j}$ and the integral over $U$ is bounded by
\[
 C\int_{\{0<R<\delta_j\}}\bigl(R^p+R^{p-1}|DG|\bigr).
\]
The integrand is integrable and the indicators tend pointwise to zero, so this error tends to zero by dominated convergence. Passing to the limit proves \eqref{eq:entropy} on $U$. On the zero set, $DG=0$ almost everywhere by \cref{lem:zeros}, which justifies the extensions.
\end{proof}

It is convenient to resolve the derivatives in directions parallel and perpendicular to $(c,s)$:
\[
 \begin{aligned}
 q_\parallel&=cq_t+sq_\theta,&q_\perp&=-sq_t+cq_\theta,\\
 \varphi_\parallel&=c\varphi_t+s\varphi_\theta,
 &\varphi_\perp&=-s\varphi_t+c\varphi_\theta.
 \end{aligned}
\]
Then \eqref{eq:amplitude} reads
\begin{equation}\label{eq:parallel}
 \varphi_\parallel=q_\perp+as,
\end{equation}
and the integrand on the right of \eqref{eq:entropy} becomes
$R^p((p-1)s\varphi_\parallel+c\varphi_\perp)$.

\begin{lemma}\label{lem:entropyenergy}
Let $p\ge2$ and suppose
\begin{equation}\label{eq:entropycondition}
 a^2+2(p+1)a-(p-1)^2-6\ge0.
\end{equation}
Then every smooth compactly supported compatible field satisfies
\[
 K_0:=\int R^{p-2}(|G_t|^2+|G_\theta|^2)\ge\int R^p=:I.
\]
\end{lemma}
\begin{proof}
Assume $I>0$, since the zero field is immediate. On $U$, $|G_j|^2=R^2(q_j^2+\varphi_j^2)$, so
\[
 K_0=\int R^p(q_t^2+q_\theta^2+\varphi_t^2+\varphi_\theta^2).
\]
This identity extends across the zero set, including at $p=2$, by \cref{lem:zeros}.

We use \eqref{eq:parallel} to split $\varphi_\parallel$ between the phase and amplitude derivatives. For a real number $\tau$, write
\[
 \varphi_\parallel=(1-\tau)\varphi_\parallel+\tau(q_\perp+as).
\]
Set
\[
 \mu=I^{-1}\int R^ps^2,\qquad
 d=pa-2-\tau a(p-1),\qquad
 e=(p-1)^2\bigl(\tau^2+(1-\tau)^2\bigr)-1.
\]
Then $0\le\mu\le1$, and \eqref{eq:entropy} gives
\[
 1+d\mu=I^{-1}\int R^p
 \bigl[(p-1)(1-\tau)s\varphi_\parallel
 +(p-1)\tau s q_\perp+c\varphi_\perp\bigr].
\]
The squared length of the coefficient vector is $1+es^2$. Weighted Cauchy--Schwarz, with the nonnegative term $q_\parallel^2$ added to the energy, yields
\begin{equation}\label{eq:entropyCS}
 |1+d\mu|^2\le\frac{K_0}{I}(1+e\mu).
\end{equation}

To make the right-hand energy at least $I$, it is enough to have $2d-e\ge0$, since
\[
 (1+d\mu)^2-(1+e\mu)=\mu(2d-e+d^2\mu).
\]
The expression $2d-e$ is a concave quadratic in $\tau$. Its maximum is attained at
\[
 \tau=\frac{p-1-a}{2(p-1)}.
\]
Substituting this choice gives
\[
 2d-e=\frac12\bigl[a^2+2(p+1)a-(p-1)^2-6\bigr],
 \qquad e=\frac{a^2+p^2-2p-1}{2}.
\]
Since $p\ge2$, $e+1/2=(a^2+p(p-2))/2\ge0$, and hence $1+e\mu\ge1/2$. Under \eqref{eq:entropycondition}, the preceding difference is nonnegative. Dividing \eqref{eq:entropyCS} by $1+e\mu$ therefore proves
\[
 \frac{K_0}{I}\ge\frac{(1+d\mu)^2}{1+e\mu}\ge1.
\]
\end{proof}

\begin{proposition}\label{prop:entropyrays}
If $p\ge2$ and
\begin{equation}\label{eq:entropyrays}
 a\le-(p+1)-\sqrt{2p^2+8}
 \quad\text{or}\quad
 a\ge-(p+1)+\sqrt{2p^2+8},
\end{equation}
then $\Lambda_p(a)=B^{p/2}$.
\end{proposition}
\begin{proof}
Let $G$ be a nonzero compactly supported compatible field and set $I=\int R^p$. On $U$,
\[
 |G_t-aG|^2=R^2((q_t-a)^2+\varphi_t^2),\qquad
 |G_\theta|^2=R^2(q_\theta^2+\varphi_\theta^2).
\]
With the probability measure $\dd\nu=R^pI^{-1}\dd t\dd\theta$, \eqref{eq:cartesian} becomes
\[
 \frac{J_a}{I}=\int\bigl((q_t-a)^2+q_\theta^2+\varphi_t^2+\varphi_\theta^2\bigr)^{p/2}\dd\nu.
\]
The power $p/2$ is at least one, so Jensen's inequality gives
\[
 \frac{J_a}{I}\ge
 \left[a^2+\frac{K_0}{I}-\frac{2a}{I}\int R^pq_t\right]^{p/2}.
\]
The last integral is zero, because $R^pq_t=p^{-1}(R^p)_t$ almost everywhere and $R^p$ is compactly supported. Condition \eqref{eq:entropyrays} is equivalent to \eqref{eq:entropycondition}. By \cref{lem:entropyenergy}, $K_0\ge I$, and therefore $J_a\ge B^{p/2}I_a$. Taking the infimum in \cref{prop:fields} and using the radial sequence proves equality. In \cref{thm:radiality}, we use the positive ray from this proposition; the calibration below gives the stated negative ray.
\end{proof}

\subsection{A homogeneous calibration}
For a scalar cylinder test, put
\[
 u=\phi_t+\beta\phi,\qquad z=\phi_\theta,\qquad \rho=(u^2+z^2)^{1/2}.
\]
Here $u$ is scalar; the vector-valued first-harmonic coefficient was used only in Section~\ref{sec:critical}. Compatibility gives
\begin{equation}\label{eq:calcompat}
 u_\theta-z=z_t-az,
\end{equation}
and
\begin{equation}\label{eq:calQ}
 Q_a\phi=(u_t-au,\sqrt2(u_\theta-z),z_\theta+u).
\end{equation}
Our aim is to find a dual field $Y$ for which $\int Q_a\phi\cdot Y=\int\rho^p$ and $|Y|\le B^{-1/2}\rho^{p-1}$. H\"older's inequality will then give the radial constant. The next identity explains how to construct such a field.

\begin{lemma}\label{lem:flux}
Suppose $H,K:\R^2\to\R$ are locally Lipschitz, continuously differentiable away from the origin, and positively homogeneous of degree $p\ge1$. Extend their gradients by zero at the origin and define
\[
 Y(u,z)=\left(H_u,\frac{H_z+K_u}{\sqrt2},K_z\right).
\]
Then
\begin{equation}\label{eq:flux}
 \int Q_a\phi\cdot Y(u,z)
 =\int\bigl[-apH(u,z)+uK_z(u,z)-zK_u(u,z)\bigr]
\end{equation}
for every $\phi\in C_c^\infty(\R\times S^1)$.
\end{lemma}
\begin{proof}
Let $U=(u,z)$. On $\{U\ne0\}$, use \eqref{eq:calcompat} in \eqref{eq:calQ} to obtain
\[
 \begin{aligned}
 Q_a\phi\cdot Y
 &=(u_t-au)H_u+(u_\theta-z)(H_z+K_u)+(z_\theta+u)K_z\\
 &=u_tH_u+z_tH_z+u_\theta K_u+z_\theta K_z
 -a(uH_u+zH_z)-zK_u+uK_z\\
 &=\partial_t(H\circ U)+\partial_\theta(K\circ U)-apH(U)-zK_u(U)+uK_z(U).
 \end{aligned}
\]
The last step uses Euler's identity for homogeneous functions. The chain rule in \cref{lem:zeros} applies almost everywhere. On $\{U=0\}$, both $DU$ and $Q_a\phi$ vanish almost everywhere, so the identity extends there. The two compositions are compactly supported, and their total derivatives integrate to zero. This proves \eqref{eq:flux}, also for $p=1$.
\end{proof}

Assume $a\ne0$. To obtain equality in the radial direction $z=0$, we choose
\begin{equation*}
 K(u,z)=\frac{\rho^{p-2}}{B}uz\qquad((u,z)\ne(0,0)).
\end{equation*}
The operator $u\partial_z-z\partial_u$ annihilates $\rho$, so
\[
 uK_z-zK_u=\frac{\rho^{p-2}}{B}(u^2-z^2).
\]
Requiring the right-hand integrand in \eqref{eq:flux} to equal $\rho^p$ now determines
\begin{equation*}
 H(u,z)=-\frac{\rho^{p-2}}{paB}
 \bigl[a^2u^2+(a^2+2)z^2\bigr].
\end{equation*}
Set $H(0,0)=K(0,0)=0$. For $p>1$, these extensions are $C^1$ at the origin. For $p=1$, their gradients are homogeneous of degree zero and bounded away from the origin. Integrating along line segments, and splitting a segment at the origin when necessary, shows that the extensions are globally Lipschitz. Thus they satisfy the hypotheses of \cref{lem:flux}.

For later use, put $T=a^2u^2+(a^2+2)z^2$. Direct differentiation away from the origin gives
\[
 \begin{aligned}
 H_u&=-\frac{\rho^{p-4}u}{paB}\bigl[(p-2)T+2a^2\rho^2\bigr],\\
 H_z&=-\frac{\rho^{p-4}z}{paB}\bigl[(p-2)T+2(a^2+2)\rho^2\bigr],\\
 K_u&=\frac{\rho^{p-4}z}{B}\bigl[\rho^2+(p-2)u^2\bigr],\\
 K_z&=\frac{\rho^{p-4}u}{B}\bigl[\rho^2+(p-2)z^2\bigr].
 \end{aligned}
\]
These formulas confirm the required identity:
\begin{equation}\label{eq:calidentity}
 -apH+uK_z-zK_u
 =\frac{\rho^{p-2}}B(T+u^2-z^2)=\rho^p.
\end{equation}
On the radial direction $z=0$, the dual field is
\[
 Y(u,0)=\frac{|u|^{p-2}u}{B}(-a,0,1),\qquad u\ne0,
\]
whose length is exactly $B^{-1/2}|u|^{p-1}$. It remains to check that the norm is no larger in other directions.

For $0\le x\le1$, define
\begin{equation}\label{eq:M}
 \begin{aligned}
 M_{p,a}(x)={}&(1-x)\left[-a^2+\left(\frac4p-2\right)x\right]^2\\
 &+\frac x2\left[-a^2+a(p-1)-\frac4p
 +\left(\frac4p-2-a(p-2)\right)x\right]^2\\
 &+a^2(1-x)[1+(p-2)x]^2.
 \end{aligned}
\end{equation}
\begin{proposition}\label{prop:calibration}
Let $p\ge1$ and $a\ne0$. If
\begin{equation}\label{eq:Mcondition}
 M_{p,a}(x)\le a^2B\qquad(0\le x\le1),
\end{equation}
then $\Lambda_p(a)=B^{p/2}$.
\end{proposition}
\begin{proof}
At a point where $\rho>0$, put $x=z^2/\rho^2$. Since $T/\rho^2=a^2+2x$, the derivative formulas give
\[
 \begin{aligned}
 \frac{aBH_u}{\rho^{p-1}}
 &=\frac u\rho\left[-a^2+\left(\frac4p-2\right)x\right],\\
 \frac{aB(H_z+K_u)}{\rho^{p-1}}
 &=\frac z\rho\left[-a^2+a(p-1)-\frac4p
 +\left(\frac4p-2-a(p-2)\right)x\right],\\
 \frac{aBK_z}{\rho^{p-1}}
 &=a\frac u\rho[1+(p-2)x].
 \end{aligned}
\]
The middle component of $Y$ has the additional factor $1/\sqrt2$. Since $u^2/\rho^2=1-x$, squaring the three components proves
\[
 \frac{a^2B^2|Y(u,z)|^2}{\rho^{2p-2}}
 =M_{p,a}\left(\frac{z^2}{\rho^2}\right).
\]
Thus \eqref{eq:Mcondition} yields $|Y|\le B^{-1/2}\rho^{p-1}$ on $\{\rho>0\}$. With the zero extension at the origin, this estimate holds everywhere for $p>1$; for $p=1$, it reads $|Y|\le B^{-1/2}$.

Let $I=\int\rho^p=\int|P_a\phi|^p$ and $J=\int|Q_a\phi|^p$. We may assume $I>0$, since the zero test is immediate. By \cref{lem:flux} and \eqref{eq:calidentity},
\[
 I=\int Q_a\phi\cdot Y.
\]
For $p>1$, H\"older's inequality with $p'=p/(p-1)$ gives
\[
 I\le J^{1/p}\left(\int|Y|^{p'}\right)^{1/p'}
 \le B^{-1/2}J^{1/p}I^{(p-1)/p}.
\]
Division and taking the $p$th power give $J\ge B^{p/2}I$. For $p=1$, the same pairing and the bound on $|Y|$ directly give $I\le B^{-1/2}J$. Taking the infimum in \eqref{eq:scalarinf} and comparing with the radial sequence prove the proposition.
\end{proof}

The norm condition is a polynomial inequality in the single directional variable $x=z^2/\rho^2$. Its verification is given in Appendix~\ref{app:polynomial}.
\begin{proposition}\label{prop:calrays}
Condition \eqref{eq:Mcondition} holds in the following regions:
\[
 \begin{aligned}
 1\le p<2:&\qquad a\le-A_0\quad\text{or}\quad a\ge A_p^+,\\
 p>2:&\qquad a\le-A_p^-.
 \end{aligned}
\]
\end{proposition}

\begin{proof}[Proof of \cref{thm:radiality}]
For $1\le p<2$, \cref{prop:calibration,prop:calrays} give $\Lambda_p(a)=B^{p/2}$ on both rays. For $p>2$, they give the negative ray, while \cref{prop:entropyrays} gives the positive ray $a\ge E_p^+$. The case $p=2$ follows from \eqref{eq:hilbert}. Taking reciprocals proves \eqref{eq:rays}.
\end{proof}

\section{Symmetry breaking and the periodic second variation}\label{sec:breaking}
We first compare a purely angular first-harmonic profile with the radial value. We then allow axial oscillation and compute the full periodic second variation. The two tests answer different questions, as the example at the end of the section shows.

\subsection{The first angular harmonic}
\begin{lemma}\label{lem:firstharmonic}
For $\Phi(\theta)=\cos\theta$ and $b=|1-a|$, put
\begin{align}
 I_1(p,a)&=\int_0^{2\pi}(b^2\cos^2\theta+\sin^2\theta)^{p/2}\dd\theta,
 \label{eq:I1}\\
 J_1(p,a)&=|a|^p\int_0^{2\pi}((b^2+1)\cos^2\theta+2\sin^2\theta)^{p/2}\dd\theta.
 \label{eq:J1}
\end{align}
Then $\Lambda_p(a)\le J_1(p,a)/I_1(p,a)$.
\end{lemma}
\begin{proof}
By \eqref{eq:angularoperators},
\[
 P_a^{\ang}\Phi=(\beta\cos\theta,-\sin\theta),\qquad
 Q_a^{\ang}\Phi=(-a\beta\cos\theta,\sqrt2a\sin\theta,-a\cos\theta).
\]
Taking squared lengths gives the two radicands in \eqref{eq:I1}--\eqref{eq:J1}. The quotient estimate follows from \cref{lem:long}.
\end{proof}

\begin{proposition}\label{prop:firstcomparison}
For every $p\ge1$ and $a\ne0$ with $-1\le a\le1/2$,
\[
 J_1(p,a)<B^{p/2}I_1(p,a).
\]
If $2<p\le4$, the same strict inequality holds for $-1-\sqrt3\le a\le-1+\sqrt3$, with $a\ne0$. For $p=2$, it is strict in the open interval and becomes equality at the endpoints.
\end{proposition}
\begin{proof}
Subtracting the squared norms in \cref{lem:firstharmonic}, we find
\[
 B|P_a^{\ang}\Phi|^2-|Q_a^{\ang}\Phi|^2
 =(1-2a)\cos^2\theta+(1-a^2)\sin^2\theta.
\]
For $-1\le a\le1/2$, the expression is nonnegative and is positive on a set of positive measure. Since $s\mapsto s^{p/2}$ is strictly increasing, integration proves the first assertion.

Now let $2<p\le4$ and write
\[
 r=\frac p2,\qquad s_0=|a|,\qquad \kappa=b^2-1.
\]
For $0\le x\le1$, define
\[
 d_p(x)=B^r(1+\kappa x)^r-s_0^{2r}(2+\kappa x)^r.
\]
Then
\[
 d_p''(x)=r(r-1)\kappa^2
 \bigl[B^r(1+\kappa x)^{r-2}-s_0^{2r}(2+\kappa x)^{r-2}\bigr]\ge0.
\]
To verify the sign, note that $a\ne1$ in the interval under consideration, so $b>0$ and
$2+\kappa x>1+\kappa x>0$. Since $r-2\le0$,
\[
 (1+\kappa x)^{r-2}\ge(2+\kappa x)^{r-2}.
\]
Also $B^r>s_0^{2r}$. The bracket is therefore positive, and $d_p$ is strictly convex when $\kappa\ne0$.

Convexity gives $d_p(x)+d_p(1-x)\ge2d_p(1/2)$. Moreover,
\[
 \begin{aligned}
 d_p(1/2)\ge0
 &\iff B(1+b^2)\ge s_0^2(3+b^2)\\
 &\iff 2-2a-a^2\ge0\\
 &\iff -1-\sqrt3\le a\le-1+\sqrt3.
 \end{aligned}
\]
In the angular integrals, $x=\cos^2\theta$, and replacing $\theta$ by $\pi/2-\theta$ replaces $x$ by $1-x$. Pairing the two angles proves $J_1\le B^{p/2}I_1$. In the open interval, $d_p(1/2)>0$, so the inequality is strict. At either endpoint, $d_p(1/2)=0$, but $\kappa\ne0$. Strict convexity then makes the paired inequality strict except where $x=1/2$, a set of measure zero. The integrated inequality is therefore still strict.

When $p=2$, one has $r=1$ and $d_2$ is affine. Thus $d_2(x)+d_2(1-x)=2d_2(1/2)$, which is positive in the open interval and zero at the endpoints.
\end{proof}

\begin{proof}[Proof of \cref{thm:breaking}]
By \cref{lem:firstharmonic,prop:firstcomparison}, the first-harmonic angular quotient is smaller than $B^{p/2}$ in each asserted strict region. The scalar cutoffs in \cref{lem:long} give admissible tests with quotient approaching that value. Hence $\Lambda_p(a)<B^{p/2}$, and taking reciprocals proves the strict comparison for $C(p,\gamma)$. At $p=2$, the endpoint equalities follow from \eqref{eq:hilbert}.
\end{proof}

\subsection{The periodic second-variation symbol}
The radial minimizing sequence has the constant field $(1,0)$ in the interior of a long plateau. We now perturb this state on a periodic cylinder. Fix $L>0$, take a smooth real function $\phi$ on $\T_L\times S^1$, and put
\[
 (V,W)=P_a\phi=(\phi_t+\beta\phi,\phi_\theta).
\]
We call $(V,W)$ an \emph{exact perturbation}, because it is obtained from a scalar potential. Define the periodic bulk deficit
\begin{equation}\label{eq:deficit}
 \begin{aligned}
 \mathcal D_{a,p,L}(\varepsilon;\phi)
 =\int_{\T_L\times S^1}\Bigl[&\bigl((\varepsilon V_t-a(1+\varepsilon V))^2
 +2\varepsilon^2(V_\theta-W)^2
 +(1+\varepsilon(V+W_\theta))^2\bigr)^{p/2}\\
 &-B^{p/2}\bigl((1+\varepsilon V)^2+\varepsilon^2W^2\bigr)^{p/2}\Bigr]\dd t\dd\theta.
 \end{aligned}
\end{equation}
At $a=1$, the base field $(1,0)$ is compatible but has no periodic scalar potential. We still use the compatible base field in \eqref{eq:deficit}, while the perturbations remain exact.

At $\varepsilon=0$, the numerator base vector is $(-a,0,1)$ and the denominator base vector is $(1,0)$. Direct differentiation gives
\[
 \mathcal D'_{a,p,L}(0;\phi)
 =pB^{p/2-1}\int_{\T_L\times S^1}(-aV_t+W_\theta)=0.
\]
Thus the second derivative is the first term that can lower the radial value. To write it, define
\begin{equation}\label{eq:coefficients}
 A=1+\frac{(p-2)a^2}{B},\qquad
 C=1-\frac{(p-2)a}{B},\qquad
 D=1+\frac{p-2}{B},\qquad
 E=(a-1)(a+2p-1).
\end{equation}
For $\xi\in(2\pi/L)\mathbb Z$ and $k\in\mathbb Z$, the Fourier symbols of $P_a$ and $Q_a$ are
\begin{align}
 p_a(\xi,k)&=(i\xi+\beta,ik),\label{eq:psymbol}\\
 q_a(\xi,k)&=\bigl((i\xi-a)(i\xi+\beta),
 \sqrt2(i\xi-a)ik,-k^2+i\xi+\beta\bigr).
 \label{eq:qsymbol}
\end{align}

\begin{lemma}\label{lem:secondsymbol}
Let $e_{\xi,k}(t,\theta)=(2\pi L)^{-1/2}e^{i(\xi t+k\theta)}$ and define the coefficient in this normalized periodic Fourier basis by
\[
 \widehat\phi(\xi,k)=\int_{\T_L\times S^1}\phi\,\overline{e_{\xi,k}}\dd t\dd\theta.
\]
Then
\begin{equation}\label{eq:secondvariation}
 \mathcal D''_{a,p,L}(0;\phi)
 =pB^{p/2-1}\sum_{\xi\in(2\pi/L)\mathbb Z}\sum_{k\in\mathbb Z}
 K_{p,a}(\xi^2,k^2)|\widehat\phi(\xi,k)|^2,
\end{equation}
where
\begin{equation}\label{eq:Ksymbol}
 K_{p,a}(x,m)=Ax^2+(A\beta^2+2Cm)x+Dm^2+Em.
\end{equation}
\end{lemma}
\begin{proof}
For $X\ne0$ and any vector $Y$,
\begin{equation}\label{eq:normvariation}
 \left.\frac{\dd^2}{\dd\varepsilon^2}|X+\varepsilon Y|^p\right|_{\varepsilon=0}
 =p|X|^{p-2}\left(|Y|^2+(p-2)\frac{(X\cdot Y)^2}{|X|^2}\right).
\end{equation}
This formula applies to both integrands in \eqref{eq:deficit}, including when $p=1$. Indeed, the periodic perturbations are bounded, so for all sufficiently small $|\varepsilon|$,
\[
 |(1+\varepsilon V,\varepsilon W)|\ge\frac12,
\]
and
\[
 \bigl|(-a,0,1)+\varepsilon(V_t-aV,\sqrt2(V_\theta-W),W_\theta+V)\bigr|
 \ge\frac{\sqrt B}{2}
\]
uniformly on the compact periodic cell. The integrands are therefore twice continuously differentiable in $\varepsilon$, with uniformly bounded derivatives, and differentiation under the integral is justified.

Apply \eqref{eq:normvariation} and expand $\phi$ in the normalized Fourier basis. At frequency $(\xi,k)$, the two amplitudes are \eqref{eq:psymbol}--\eqref{eq:qsymbol}. After removing the common factor $pB^{p/2-1}|\widehat\phi(\xi,k)|^2$, the contribution is
\[
 |q_a|^2+\frac{p-2}{B}|-a(q_a)_1+(q_a)_3|^2
 -B\bigl(|p_a|^2+(p-2)|(p_a)_1|^2\bigr).
\]
To collect the terms, put $x=\xi^2$ and $m=k^2$. The symbols give
\[
 \begin{aligned}
 |p_a|^2&=x+\beta^2+m,\quad |(p_a)_1|^2=x+\beta^2,\\
 |q_a|^2&=(x+a\beta)^2+x(\beta-a)^2
 +2m(x+a^2)+x+(\beta-m)^2,\\
 |-a(q_a)_1+(q_a)_3|^2
 &=(ax+\beta B-m)^2+x[1-a(\beta-a)]^2.
 \end{aligned}
\]
Substitution shows that the coefficients of $x^2,xm,m^2,x,m$ are $A,2C,D,A\beta^2,E$, respectively; the constant term is zero. The quadratic form has constant real coefficients. Orthogonality removes cross terms between distinct Fourier modes, and Parseval's identity gives \eqref{eq:secondvariation}.
\end{proof}

For the first angular harmonic, define
\[
 F_p(a)=D+E,\qquad L_p(a)=A\beta^2+2C,
\]
and, writing $r_- =\max\{-r,0\}$ for a real number $r$, put
\begin{equation}\label{eq:Sp}
 S_p(a)=F_p(a)-\frac{(L_p(a)_-)^2}{4A}.
\end{equation}
The quantity $A=(1+(p-1)a^2)/B$ is positive, so this expression is defined for every $p\ge1$ and $a\in\R$.

\begin{lemma}\label{lem:firstcontrols}
The symbol in \eqref{eq:Ksymbol} satisfies $K_{p,a}(x,k^2)\ge0$ for all $x\ge0$ and $k\in\mathbb Z$ if and only if $K_{p,a}(x,1)\ge0$ for every $x\ge0$. Moreover,
\[
 \min_{x\ge0}K_{p,a}(x,1)=S_p(a).
\]
\end{lemma}
\begin{proof}
For $k=0$, one has $K_{p,a}(x,0)=Ax(x+\beta^2)\ge0$. Now let $m=k^2\ge1$ and write $x=my$. Then
\begin{equation}\label{eq:modereduction}
 \frac{K_{p,a}(my,m)}m
 =m(Ay^2+2Cy+D)+A\beta^2y+E.
\end{equation}
We first check that $Ay^2+2Cy+D\ge0$ for $y\ge0$. Besides $A>0$, we have $D=(a^2+p-1)/B\ge0$. If $1\le p\le2$, then
\[
 C=1+(2-p)\frac aB\ge1-\frac{2-p}{2}>0,
\]
because $a/B\ge-1/2$. For $p\ge2$, the sign is immediate when $C\ge0$. When $C<0$, the identity
\[
 AD-C^2=\frac{(p-2)(a+1)^2}{B}\ge0
\]
shows that the quadratic is nonnegative on all of $\R$.

It follows from \eqref{eq:modereduction} that
\[
 \frac{K_{p,a}(my,m)}m
 =K_{p,a}(y,1)+(m-1)(Ay^2+2Cy+D).
\]
Thus nonnegativity for the first angular harmonic implies it for all $m=k^2\ge1$. The converse follows by taking $k=1$.

Finally, $K_{p,a}(x,1)=Ax^2+L_p(a)x+F_p(a)$. If $L_p(a)\ge0$, its minimum on $[0,\infty)$ is $F_p(a)$. Otherwise, the minimum is at $x=-L_p(a)/(2A)$ and equals $F_p(a)-L_p(a)^2/(4A)$. Both cases are expressed by \eqref{eq:Sp}.
\end{proof}

\begin{theorem}\label{thm:variation}
Let $p\ge1$ and $a\in\R$. The periodic bulk second variation at $(1,0)$ is nonnegative against every smooth exact perturbation on every periodic cylinder if and only if $S_p(a)\ge0$. If $S_p(a)<0$, then
\[
 C(p,\gamma)>B^{-p/2}.
\]
For $1\le p\le4$, one has $S_p(a)=F_p(a)$ and
\begin{equation}\label{eq:BF}
 BF_p(a)=a(a-1)(a^2-a+2)+p(2a^3-2a^2+2a-1).
\end{equation}
\end{theorem}
\begin{proof}
For a fixed period $L$, \cref{lem:secondsymbol} shows that nonnegativity for every exact perturbation is equivalent to
\[
 K_{p,a}(\xi^2,k^2)\ge0\qquad
 \left(\xi\in\frac{2\pi}L\mathbb Z,\ k\in\mathbb Z\right).
\]
As $L$ varies, every $x>0$ is represented by $\xi^2$, by choosing $L=2\pi/\sqrt x$ and $\xi=2\pi/L$. The case $x=0$ corresponds to zero axial frequency. The first assertion therefore follows from \cref{lem:firstcontrols}.

Suppose $S_p(a)<0$, and let $x_*$ minimize $K_{p,a}(x,1)$ on $[0,\infty)$. If $x_*>0$, choose $L=2\pi/\sqrt{x_*}$ and $\xi=\sqrt{x_*}$. If $x_*=0$, choose any $L$ and put $\xi=0$. Then
\[
 \psi(t,\theta)=\cos(\xi t+\theta)
\]
is a real periodic potential with $\mathcal D''_{a,p,L}(0;\psi)<0$.

To turn this perturbation into a periodic scalar test about the radial state, we need $a\ne1$. This follows from the criterion itself. At $a=1$, one has $\beta=0$, $F_p(1)=p/2$, and $L_p(1)=4-p$, whence
\[
 S_p(1)=
 \begin{cases}
 p/2,&1\le p\le4,\\
 4-8/p,&p>4.
 \end{cases}
\]
Both expressions are positive. Thus $S_p(a)<0$ implies $\beta\ne0$, and the constant potential $\phi_0=1/\beta$ satisfies $P_a\phi_0=(1,0)$.

Consider on the periodic cell the quotient
\[
 \mathcal R(\varepsilon)
 =\frac{\int|Q_a(\phi_0+\varepsilon\psi)|^p}
 {\int|P_a(\phi_0+\varepsilon\psi)|^p}.
\]
The uniform separation from zero used in \cref{lem:secondsymbol} makes its numerator $N(\varepsilon)$ and denominator $D(\varepsilon)$ twice continuously differentiable near zero, including at $p=1$. Their deficit is $N(\varepsilon)-B^{p/2}D(\varepsilon)=\mathcal D_{a,p,L}(\varepsilon;\psi)$. It has zero value and first derivative at zero, and negative second derivative. Since $D(0)=2\pi L>0$, differentiation of the quotient gives
\[
 \mathcal R(0)=B^{p/2},\qquad \mathcal R'(0)=0,\qquad
 \mathcal R''(0)=\frac{\mathcal D''_{a,p,L}(0;\psi)}{D(0)}<0.
\]
Hence $\mathcal R(\varepsilon)<B^{p/2}$ for all sufficiently small nonzero $\varepsilon$. Applying \cref{lem:periodic} to this fixed scalar periodic profile gives compactly supported tests with quotient below $B^{p/2}$. Taking reciprocals proves strict symmetry breaking.

For $1\le p\le4$,
\[
 C=\frac{a^2-(p-2)a+1}{B}\ge0,
\]
because the discriminant of the numerator is $(p-2)^2-4\le0$. Thus $L_p(a)=A\beta^2+2C\ge0$, and \eqref{eq:Sp} gives $S_p(a)=F_p(a)$. Expanding $B(D+E)$ proves \eqref{eq:BF}.
\end{proof}

For $p>4$, the axial correction in \eqref{eq:Sp} can change the sign. In \cref{prop:axial}, we verify that $F_6(7/10)>0>S_6(7/10)$. Thus zero axial frequency is not enough for the full second-variation test.

The next consequence makes the distinction between the quadratic test and global symmetry breaking explicit.
\begin{corollary}\label{cor:nonnecessity}
At $(p,a)=(4,7/10)$, the periodic bulk second variation is nonnegative for every smooth exact perturbation and every period, but
\[
 C(4,6/5)>\left(\frac{100}{149}\right)^2=C_{\rad}(4,6/5).
\]
In particular, a negative periodic second variation is not necessary for strict symmetry breaking.
\end{corollary}
\begin{proof}
Substituting $p=4$ and $a=7/10$ in \eqref{eq:BF}, and using $B=149/100$, gives
\[
 S_4(7/10)=F_4(7/10)=\frac{481}{14900}>0.
\]
The second-variation assertion follows from \cref{thm:variation}. On the other hand,
$-1-\sqrt3<7/10<-1+\sqrt3$, since $(17/10)^2<3$. Thus \cref{thm:breaking}(ii) gives the strict comparison with $B^{-2}$. Finally, $\gamma=a+1-2/p=6/5$.
\end{proof}

Nonnegative second variation alone does not decide nonlinear local minimality. Here the strict global comparison comes from an angular profile of finite size, rather than a negative quadratic perturbation at the radial state.

\section{A nonlinear interaction between angular modes}\label{sec:interaction}
The second variation separates into Fourier modes, but the nonlinear angular quotient does not. We illustrate this at $p=4$ and $a=-1$, where a combination of the first and third harmonics improves every nonzero real single-frequency angular profile.

Here $\beta=2$ and $B=2$. For a nonzero real $\Phi\in C^\infty(S^1)$, define
\[
 N[\Phi]=\int_0^{2\pi}\bigl(4\Phi^2+(\Phi')^2\bigr)^2\dd\theta,
 \qquad
 D[\Phi]=\int_0^{2\pi}\bigl(4\Phi^2+2(\Phi')^2+(\Phi''+2\Phi)^2\bigr)^2\dd\theta.
\]
By \cref{lem:long}, $N[\Phi]/D[\Phi]\le C(4,-1/2)$. Since
\[
 \cos^3\theta=\frac{3\cos\theta+\cos3\theta}{4},
\]
the third harmonic is a natural first choice for an interaction with $\cos\theta$.

\begin{proposition}\label{prop:interaction}
Let $\Phi_\varepsilon(\theta)=\cos\theta+\varepsilon\cos3\theta$. Then
\begin{align}
 \frac{N[\Phi_\varepsilon]}{2\pi}
 &=\frac{59}{8}+\frac{21}{2}\varepsilon+\frac{115}{2}\varepsilon^2
 +\frac{363}{8}\varepsilon^4,\label{eq:Npoly}\\
 \frac{D[\Phi_\varepsilon]}{2\pi}
 &=\frac{107}{8}+\frac92\varepsilon+\frac{587}{2}\varepsilon^2
 +\frac{11307}{8}\varepsilon^4.\label{eq:Dpoly}
\end{align}
The derivative of $N[\Phi_\varepsilon]/D[\Phi_\varepsilon]$ at $\varepsilon=0$ is positive, and
\[
 \frac{N[\Phi_{1/20}]}{D[\Phi_{1/20}]}>\frac{14}{25}.
\]
Every nonzero real single-frequency angular profile has quotient at most $59/107$, with equality for the first harmonic.
\end{proposition}
\begin{proof}
Write $C_j=\cos(j\theta)$. Product formulas give
\[
 \begin{aligned}
 4\Phi_\varepsilon^2+(\Phi_\varepsilon')^2
 &=\frac{5+13\varepsilon^2}{2}
 +\left(\frac32+7\varepsilon\right)C_2+\varepsilon C_4
 -\frac52\varepsilon^2C_6,\\
 4\Phi_\varepsilon^2+2(\Phi_\varepsilon')^2+(\Phi_\varepsilon''+2\Phi_\varepsilon)^2
 &=\frac{7+71\varepsilon^2}{2}
 +\left(\frac32+3\varepsilon\right)C_2-9\varepsilon C_4
 +\frac{35}{2}\varepsilon^2C_6.
 \end{aligned}
\]
Each $C_j$ has mean zero, and $\int_0^{2\pi}C_jC_\ell\dd\theta=\pi\delta_{j\ell}$ for $j,\ell\ge1$. Squaring and integrating therefore proves \eqref{eq:Npoly}--\eqref{eq:Dpoly}. In particular,
\[
 \left.\frac{\dd}{\dd\varepsilon}\frac{N[\Phi_\varepsilon]}{D[\Phi_\varepsilon]}\right|_{\varepsilon=0}
 =\frac{6864}{11449}>0.
\]
Since the quotient at zero is $59/107$, it is strictly larger for all sufficiently small positive $\varepsilon$. At $\varepsilon=1/20$, the exact calculation
\[
 \frac{25N[\Phi_{1/20}]-14D[\Phi_{1/20}]}{2\pi}
 =\frac{389977}{1280000}>0
\]
proves the stated lower bound.

It remains to compare with the other single-frequency profiles. Every nonzero real profile supported at frequency $k\ge1$ has the form $R_0\cos(k\theta-\delta)$. Scaling and angular translation leave the quotient unchanged, so it suffices to take $\Phi_k=\cos(k\theta)$. With $x=k^2$,
\[
 4\Phi_k^2+(\Phi_k')^2=x+(4-x)\cos^2(k\theta)
\]
and
\[
 4\Phi_k^2+2(\Phi_k')^2+(\Phi_k''+2\Phi_k)^2
 =2x+(x^2-6x+8)\cos^2(k\theta).
\]
Using $\int_0^{2\pi}\cos^2(k\theta)\dd\theta=\pi$ and
$\int_0^{2\pi}\cos^4(k\theta)\dd\theta=3\pi/4$, we obtain
\[
 \frac{N[\Phi_k]}{D[\Phi_k]}
 =\frac{3x^2+8x+48}{3x^4-20x^3+92x^2-160x+192}.
\]
At $k=1$, the value is $59/107$. For $k\ge2$, $x\ge4$, and
\begin{multline*}
 59(3x^4-20x^3+92x^2-160x+192)-107(3x^2+8x+48)\\
 =(x-1)(177x^3-1003x^2+4104x-6192)>0.
\end{multline*}
Indeed, the derivative of the cubic is $531x^2-2006x+4104$. Its discriminant is negative and its leading coefficient is positive, so the cubic is increasing on $\R$. Its value at $x=4$ is $5504$, which proves the sign on $[4,\infty)$. The constant harmonic has quotient $1/4$. This completes the comparison.
\end{proof}

Since $14/25-59/107=23/2675>0$, the explicit profile $\Phi_{1/20}$ improves every nonzero real single-frequency angular profile. In particular, $C(4,-1/2)>14/25$. Thus optimizing the angular quotient one Fourier frequency at a time would miss even this elementary two-frequency improvement.

\section{The remaining problem}\label{sec:remaining}
The critical law determines the leading behaviour of the unrestricted constant as $a\to0$, and \cref{thm:radiality} determines its exact value on exterior rays. For each fixed $p\ne2$, the parameters not settled by these results lie in
\[
 (-A_0,A_p^+)\quad\text{if }1\le p<2,
 \qquad
 (-A_p^-,E_p^+)\quad\text{if }p>2,
\]
with the value at $a=0$ already known to be infinite. The exterior endpoints arise from sufficient entropy and calibration conditions. Determining the exact transition points remains part of the problem.

Within a symmetry-breaking region, one must still determine the unrestricted nonlinear infimum. The quadratic criterion does not settle this question: \cref{cor:nonnecessity} shows that global symmetry breaking can occur even when every periodic second variation is nonnegative. Purely angular tests provide a different sufficient criterion, namely
\[
 \int_0^{2\pi}|Q_a^{\ang}\Phi|^p\dd\theta
 <B^{p/2}\int_0^{2\pi}|P_a^{\ang}\Phi|^p\dd\theta,
 \qquad P_a^{\ang}\Phi\not\equiv0.
\]
Whether optimization over all such angular profiles determines the unrestricted infimum in the transition region is not answered here.

The examples clarify two separate limitations. \Cref{prop:interaction} shows that even the purely angular nonlinear problem cannot be optimized one frequency at a time. \Cref{prop:axial} shows that axial frequency cannot be discarded from the linearized problem at the radial state. The latter does not prove that axial oscillation is necessary to obtain the full nonlinear sharp constant. That comparison, together with the exact value throughout the remaining transition interval, is left open.

\appendix
\section{Algebraic verifications}\label{app:algebra}
We give the polynomial signs needed for the calibration rays and then verify the axial-frequency example used in Section~\ref{sec:breaking}.

\subsection{The polynomial test for the calibration rays}\label{app:polynomial}
\begin{lemma}\label{lem:polynomial}
For $M_{p,a}$ in \eqref{eq:M},
\begin{equation}\label{eq:qfactor}
 M_{p,a}(x)-a^2B=xq_{p,a}(x),\qquad q_{p,a}(x)=m_0+m_1x+m_2x^2,
\end{equation}
where
\begin{equation}\label{eq:mcoeffs}
 m_0=-\frac{P_p(a)}{2p^2},\qquad
 m_1=\frac{(p-2)(ap-2)(a^2-3a-2)}p,\qquad
 m_2=-\frac{(p-2)^2(ap-2)^2}{2p^2}\le0.
\end{equation}
Moreover,
\begin{align}
 q_{p,a}(1)&=-\frac12(a-1)(a^3+3a^2+4),\label{eq:qendpoint}\\
 q_{p,a}'(1)&=-\frac{(p-2)(ap-2)}{p^2}T_p(a),
 \qquad T_p(a)=-pa^2+ap^2+ap+4.
 \label{eq:qderivative}
\end{align}
\end{lemma}
\begin{proof}
Expanding \eqref{eq:M} and collecting powers of $x$ gives
\[
 \begin{aligned}
 2p^2(M_{p,a}(x)-a^2B)
 =x\bigl[&-P_p(a)+2p(p-2)(ap-2)(a^2-3a-2)x\\
 &-(p-2)^2(ap-2)^2x^2\bigr].
 \end{aligned}
\]
Dividing by $2p^2$ proves \eqref{eq:qfactor}--\eqref{eq:mcoeffs}. Substitution of $x=1$ gives \eqref{eq:qendpoint}, and differentiating gives \eqref{eq:qderivative}.
\end{proof}

\begin{proof}[Proof of \cref{prop:calrays}]
Write $q=q_{p,a}$. By \eqref{eq:mcoeffs}, it is a concave quadratic. We check its sign separately on the three stated rays.

Suppose first that $1\le p<2$ and $a\le-A_0$. Write $a=-A$. The factorization
\[
 A^4-2A^3-3A^2-4A-4=(A+1)(A^3-3A^2-4)
\]
reduces \eqref{eq:A0} to $g(A)=A^3-3A^2-4=0$ for $A>1$. Since $g'(A)=3A(A-2)$, $g(1)=-6$, and $g(2)=-8$, the function is negative on $[1,2]$. It is strictly increasing on $(2,\infty)$ and tends to infinity. Thus there is exactly one root $A_0>2$, as claimed in the introduction. For $A\ge A_0$,
\[
 a^3+3a^2+4=-A^3+3A^2+4\le0.
\]
Since $a-1<0$, \eqref{eq:qendpoint} gives $q(1)\le0$. Also,
\[
 \partial_pT_p(a)=-a^2+a(2p+1)<0\qquad(a<0),
\]
so on this ray
\[
 T_p(a)\le T_1(a)=-a^2+2a+4<0.
\]
Both $p-2$ and $ap-2$ are negative. Hence \eqref{eq:qderivative} gives $q'(1)\ge0$. Concavity implies $q'(x)\ge q'(1)\ge0$ for $0\le x\le1$, and therefore $q(x)\le q(1)\le0$.

Next suppose $p>2$ and $a\le-A_p^-$. The polynomial in \eqref{eq:thresholds} satisfies
\[
 P_p(-1)=-p^4-4p^3-4(p-2)^2<0,
\]
whereas $P_p(-A)\to\infty$ as $A\to\infty$. Its largest positive root $A_p^-$ therefore exists and is greater than one. The polynomial has positive leading coefficient and no zero beyond this root, so $P_p(-A)\ge0$ for $A\ge A_p^-$. Thus $m_0\le0$. Also $p-2>0$, $ap-2<0$, and $a^2-3a-2>0$ for $a<-1$. Hence $m_1\le0$. With $m_2\le0$, the quadratic is nonincreasing on $[0,1]$, and $q(x)\le q(0)=m_0\le0$.

Finally, let $1\le p<2$ and $a\ge A_p^+$. Since
\[
 P_p(1)=-(p-2)^2(p+2)^2<0
\]
and $P_p(A)\to\infty$ as $A\to\infty$, the largest positive root $A_p^+$ exists and exceeds one. As above, $P_p(a)\ge0$ for $a\ge A_p^+$, so $m_0\le0$. Equation \eqref{eq:qendpoint} also gives $q(1)\le0$, since $a>1$.

If $m_1\le0$, the quadratic is nonincreasing and $q\le q(0)\le0$. If $q'(1)\ge0$, it is nondecreasing and $q\le q(1)\le0$. We are left with
\begin{equation}\label{eq:remainingcase}
 m_1>0,\qquad q'(1)<0.
\end{equation}
These inequalities force
\[
 ap-2>0,\qquad a^2-3a-2<0,\qquad T_p(a)<0.
\]
Indeed, $ap-2=0$ would make $m_1=0$. If $ap-2<0$, then $a<2/p\le2$, whereas $m_1>0$ would force $a^2-3a-2>0$, and hence $a>(3+\sqrt{17})/2>2$, a contradiction. Therefore $ap-2>0$. Since $p-2<0$, the first inequality in \eqref{eq:remainingcase} now implies $a^2-3a-2<0$, and the second one, together with \eqref{eq:qderivative}, implies $T_p(a)<0$.

The discriminant of $q$ is
\begin{equation}\label{eq:discriminant}
 m_1^2-4m_2m_0
 =\frac{(p-2)^2(ap-2)^2}{p^4}R_p(a),
\end{equation}
where the following expansion is useful for its sign:
\[
 \begin{aligned}
 R_p(a)={}&T_p(a)(2ap^2+4ap+p^3+4p^2+8)\\
 &-a(p^5+5p^4+4p^3+12p^2+16p)-4p^3-12p^2-16.
 \end{aligned}
\]
For $a>1$ and $p\ge1$, the factor multiplying $T_p(a)$ and the polynomial multiplying $a$ are positive. Since $T_p(a)<0$, all three displayed contributions to $R_p(a)$ are negative. Thus \eqref{eq:discriminant} is negative. Moreover, $m_2<0$ in the case under consideration, so $q$ is strictly negative on $\R$. This proves $q\le0$ on $[0,1]$ in the last case as well. By \eqref{eq:qfactor}, the desired condition \eqref{eq:Mcondition} holds on all three rays.
\end{proof}

\subsection{An axial-frequency example}
\begin{proposition}\label{prop:axial}
At $(p,a)=(6,7/10)$,
\[
 F_6(7/10)=\frac{2601}{14900}>0,
 \qquad S_6(7/10)=-\frac{6977641}{82248000}<0.
\]
The first-mode symbol $x\mapsto K_{6,7/10}(x,1)$ attains its minimum at
\[
 x=\xi^2=\frac{4619}{13800}>0.
\]
All modes with zero axial frequency have nonnegative second variation.
\end{proposition}
\begin{proof}
Substitution in \eqref{eq:coefficients} gives
\[
 A=\frac{345}{149},\qquad C=-\frac{131}{149},\qquad
 D=\frac{549}{149},\qquad
 L_6(7/10)=-\frac{31}{20},\qquad F_6(7/10)=\frac{2601}{14900}.
\]
Since $L_6(7/10)<0$, the quadratic is minimized at
\[
 x_*=-\frac{L_6(7/10)}{2A}=\frac{4619}{13800}.
\]
Formula \eqref{eq:Sp} gives
\[
 S_6(7/10)=F_6(7/10)-\frac{L_6(7/10)^2}{4A}
 =-\frac{6977641}{82248000}.
\]
On the other hand,
\[
 K_{6,7/10}(0,m)=m\bigl[F_6(7/10)+D(m-1)\bigr]>0\qquad(m\ge1),
\]
and the constant mode has symbol zero. Thus the zero-axial-frequency modes are all nonnegative, while $K_{6,7/10}(x_*,1)<0$. Taking $L=2\pi/\sqrt{x_*}$ makes $\xi=\sqrt{x_*}$ an admissible periodic frequency and gives a negative mode.
\end{proof}

\section*{Acknowledgement} The author thanks Jian-Feng Zhu for useful discussions during the preparation of the manuscript.

\section*{AI Declaration}
ChatGPT (OpenAI) was used for language editing, organisation, and assistance with finite symbolic calculations. The author independently checked all mathematical contents and takes full responsibility for the final manuscript.

\end{document}